\documentclass[a4paper, 12pt]{amsart}

\usepackage{amsfonts}
\usepackage{amssymb}
\usepackage{latexsym}
\usepackage{graphics}
\usepackage[2emode]{psfrag}
\usepackage{mathrsfs}
\usepackage{amsthm}
\usepackage{amsmath}
\usepackage[all]{xy}
\usepackage{enumerate}
\usepackage[latin1]{inputenc}
\usepackage{lscape}
\usepackage{a4wide}
\usepackage[usenames]{color}
\usepackage{pdfsync}

\input{xy}
 \xyoption{all}
\newtheorem{theorem}{\sc Theorem}[section]
\newtheorem{proposition}[theorem]{\sc Proposition}

\newtheorem{lemma}[theorem]{\sc Lemma}
\newtheorem{corollary}[theorem]{\sc Corollary}
\theoremstyle{definition}

\newtheorem{example}[theorem]{\sc Example}

\theoremstyle{remark}
\newtheorem{remark}[theorem]{\sc Remark}

\newcommand{\tensor}[1]{\otimes_{#1}}

\newcommand{\env}[1]{#1\tensor{}#1^{o}}
\newcommand{\Sf}[1]{\mathsf{#1}}
\newcommand{\fk}[1]{\mathfrak{#1}}
\newcommand{\cat}[1]{\mathcal{#1}}
\newcommand{\rmod}[1]{\Sf{Mod}_{#1}}
\newcommand{\lmod}[1]{{}_{#1}\Sf{Mod}}
\newcommand{\bimod}[1]{{}_{#1}{\Sf{ Mod}}_{#1}}

\renewcommand{\hom}[3]{\mathrm{Hom}_{#1}\left(#2,\,#3\right)}

\newcommand{\td}[1]{\widetilde{#1}}

\renewcommand{\Re}{R^{\Sf{e}}}
\newcommand{\bara}[1]{\overline{#1}}

\newcommand{\op}[1]{#1^o}
\newcommand{\lr}[1]{\left(\underset{}{} #1 \right)}
\newcommand{\Lr}[1]{\left[\underset{}{} #1 \right]}

\newcommand{\cross}[2]{#1 \times_{R}\, #2}
\newcommand{\ring}[1]{#1\text{-}\mathsf{Rings}}
\newcommand{\End}[2]{\mathrm{End}_{#1}(#2)}

\newcommand{\coring}[1]{\mathfrak{#1}}
\newcommand{\lcomod}[1]{ {}_{#1}\mathsf{Comod}}
\newcommand{\can}[1]{\mathsf{can}_{#1}}
\newcommand{\chain}[1]{\mathit{Ch}_{+}(#1)}

\newcommand{\dD}{\mathscr{D}}

\newcommand{\gG}{\mathscr{G}}

\newcommand{\iI}{\mathscr{I}}
\newcommand{\jJ}{\mathscr{J}}

\newcommand{\lL}{\mathscr{L}}

\newcommand{\oO}{\mathscr{O}}

\newcommand{\tT}{\mathscr{T}}

\begin{document}
\allowdisplaybreaks

\title{Categories of comodules and  chain complexes of modules}
\author{A. Ardizzoni}
\address{University of Ferrara, Department of Mathematics, Via Machiavelli
35, Ferrara, I-44121, Italy}
\email{rdzlsn@unife.it}
\urladdr{http://www.unife.it/utenti/alessandro.ardizzoni}
\author{L. El Kaoutit}
\address{Universidad de Granada, Departamento de \'{A}lgebra. Facultad de Educaci\'{o}n y Humanidades de Ceuta.
El Greco N. 10. E-51002 Ceuta, Spain}
\email{kaoutit@ugr.es}
\author{C. Menini}
\address{University of Ferrara, Department of Mathematics, Via Machiavelli
35, Ferrara, I-44121, Italy}
\email{men@unife.it}
\urladdr{http://www.unife.it/utenti/claudia.menini}
\date{\today}
\thanks{This paper was written while the first and the last authors were members of
G.N.S.A.G.A. with partial financial support from M.I.U.R. within
the National Research Project PRIN 2007. The second author was
supported by grant MTM2010-20940-C02-01 from the Ministerio de
Educaci\'{o}n y Ciencia of Spain.  His stay, as a visiting
professor at University of Ferrara, was supported by I.N.D.A.M..}

\begin{abstract}
Let $\lL(A)$ denote the  coendomorphism left $R$-bialgebroid associated to a left  finitely generated and projective extension of rings  $R \to A$ with identities. We show that the category of left comodules over an epimorphic image of $\lL(A)$ is equivalent to the    category of chain complexes of left $R$-modules. This equivalence is monoidal whenever
$R$ is commutative and $A$ is an $R$-algebra. This is a generalization, using entirely new tools, of results by  B. Pareigis and D. Tambara for chain complexes of vector spaces over fields. Our approach relies heavily on the non commutative theory of
Tannaka reconstruction, and the generalized faithfully flat
descent for small additive categories, or rings with enough
orthogonal idempotents.
\end{abstract}

\maketitle


\pagestyle{headings}

\section{Introduction}

\subsection{Methodology and motivation overviews.}
The starting point of this paper is a result due to B. Pareigis
\cite[Theorem 18]{Pareigis:81} which  asserts that the category of
unbounded complexes of vector spaces is monoidally equivalent to
the category of left comodules over a certain Hopf algebra which
is neither commutative nor cocommutative. Later on, in
\cite[Theorem 4.4]{Tambara}, D. Tambara associated to every finite
dimensional algebra $A$ over a field $\Bbbk$, a bialgebra
$\lL(A)$ (termed \emph{coendomorphism bialgebra}) such that the
category of left comodules $\lcomod{\lL(A)}$ is monoidally
equivalent to the category $\chain{\Bbbk}$ of
chain complexes of $\Bbbk$-vector spaces. The Hopf algebra
considered by B. Pareigis is recovered by choosing $A=\Bbbk\oplus
\Bbbk t$ with $t^2=0$, i.e. the trivial extension of $\Bbbk$, and
localizing the bialgebra $\lL(A)$ using a multiplicative set
generated by a single grouplike element. The equivalence of
categories established by Tambara relies on the use of a
variant of the equivalence between simplicial $\Bbbk$-vector
spaces and chain complexes of $\Bbbk$-vector spaces, provided by
the normalization functor, due to Dold and Kan, see \cite[Theorem
1.9, Corollary 1.12]{Dold:1958} or \cite[Theorem
2.4]{Goerss-Jardine}. The functor that provides such equivalence
is given, in some sense, by tensoring chain complexes with the
augmented cochain complex $Q_{\bullet}$ constructed using the
Amitsur cosimplicial vector space attached to the $\Bbbk$-algebra
$A$. Note that $Q_{\bullet}$ is the universal differential graded $\Bbbk$-algebra of $A$, given by $Q_0=\Bbbk$, $Q_1=A$ and
$Q_n=K\tensor{A} \cdots \tensor{A} K$, $(n-1)$-times for $n \geq
2$, where $K$ is the kernel of the multiplication of $A$. The
construction of this functor will be clarified in Section
\ref{Main Results}, see also the forthcoming  subsection. A
different approach to Pareigis's result, using Tannaka
reconstruction for several-objects coalgebras, was also given by
Paddy McCrudden in \cite[Examples 6.6, 6.9]{McCrudden:2000}, where
the same coendomorphism bialgebra $\lL(A)$ was constructed for
a commutative base ring $\Bbbk$ instate of a field.

Let $R$ be an algebra over a commutative ring $\Bbbk$. The purpose
of this paper is to investigate 
the
relationship between the category  of left comodules over certain left
$R$-bialgebroids, termed \emph{coendomorphism bialgebroids} coming from the $\times_R$-bialgebra defined in \cite[Remark 1.7]{Tambara}, and the category of chain complexes of left
$R$-modules. Tambara's results, and in
particular Pareigis's one, are then immediate consequences of the
general theory developed here. It is noteworthy that our methods can
be seen as new and more conceptual even for the case of vector
spaces. Indeed, we will see  why concretely the trivial extension of rings, already  considered by  B. Pareigis, induces the above equivalence of categories.
Our approach makes use of the ''non commutative'' Tannakian
categories theory  following  the spirit of \cite{Deligne:1990,
Bruguieres:1994} and \cite{P.H. Hai:2008}, as well as of the
generalized faithfully flat descent for rings with enough
orthogonal idempotents stated in \cite{El Kaoutit/Gomez:2004}. We
mean that all  (left) bialgebroids arising here come in fact from
the non commutative version of Tannaka reconstruction process
which in our approach involves rings with enough orthogonal idempotents.

In the setting of  non commutative Tannakian categories, one basically starts with a small
$\Bbbk$-linear monoidal category $(\cat{A}, \tensor{}, {\bf 1})$
and a faithful monoidal functor \footnote[3]{Our setting requires
an isomorphism only at the level of unit. That is, $R\cong \omega
({\bf 1})$, while $\omega(-\tensor{}-) \to
\omega(-)\tensor{R}\omega(-)$ is  not necessarily a natural
isomorphism.} from $\cat{A}$ to the category of $R$-bimodules,
$\omega: \cat{A} \to \bimod{R}$ (the fiber functor), valued in the
category $add({}_RR)$ of finitely generated and projective left $R$-modules (i.e.
locally free sheaves of finite rank). There are several objects
under consideration:
$$\Sigma(\omega)= \underset{\fk{p} \in \cat{A}}{\oplus} \omega(\fk{p}),\quad
{}^{\vee}\Sigma(\omega)= \underset{\fk{p} \in \cat{A}}{\oplus}
{}^*\omega(\fk{p}),\qquad \gG(\cat{A})=\underset{\fk{p},\, \fk{p}'
\in \cat{A}}{\oplus}\hom{\cat{A}^{o}}{\fk{p}}{\fk{p}'}. $$ Here
the second is the right $R$-module direct sum of the left duals
while the third is Gabriel's ring with enough orthogonal
idempotents, introduced in \cite{Gabriel:1962}, attached to the
opposite category $\cat{A}^o$ of $\cat{A}$. Using the canonical actions, we consider
$\cat{L}(\omega) := \Sigma(\omega)\tensor{\gG(A)}
{}^{\vee}\Sigma(\omega)$ as an $\Re$-bimodule, where
$\Re:=R\tensor{\Bbbk}R^o$ denotes the enveloping ring.
A well known argument in small additive
categories  says that the object $\cat{L}(\omega)$ solves the
following universal problems in $R$-bimodules
\begin{eqnarray*}
{\rm Nat}\lr{\omega,\, -\tensor{R}\omega} &\cong & \hom{R-R}{\cat{L}(\omega)}{-}, \\
{\rm Nat}\lr{\omega\tensor{R}\omega,\, -\tensor{R}(\omega\tensor{R}\omega)} & \cong & \hom{R-R}{\cat{L}(\omega)\tensor{\Re}\cat{L}(\omega)}{-},
\end{eqnarray*}
where the $R$-bimodule structures of  $\cat{L}(\omega)$ have been
chosen properly. It is indeed this solution which allows us to
construct a left $R$-bialgebroid (or a Hopf bialgebroid if
desired). Of course there is an obvious (monoidal) functor
connecting left unital $\gG(\cat{A})$-modules and left
$\cat{L}(\omega)$-comodule, namely
$$\Sigma(\omega)\tensor{\gG(\cat{A})}- : \lmod{\gG(\cat{A})} \longrightarrow \lcomod{\cat{L}(\omega)}.$$
In the case when each of the left $R$-modules $\omega(\fk{p})$ is
endowed with a structure of left $\coring{C}$-comodule for some
$R$-coring $\coring{C}$ (or certain left $R$-bialgebroid), there
is a map of $R$-corings, known as a \emph{canonical map},
$$ \can{\gG(\cat{A})} : \cat{L}(\omega) \longrightarrow \coring{C}$$
defined by using the left $\coring{C}$-coaction of the
$\omega(\fk{p})$'s.  This homomorphism of corings is not in general bijective, see \cite{El
Kaoutit/Gomez:2004} for more discussions. The associated coinduction functor of the canonical map leads to
the following composition of functors
\begin{equation*}\label{Eq-0}
\xymatrix@C=40pt{ \lmod{\gG(\cat{A})} \ar@{->}^-{\Sigma(\omega)\tensor{\gG(\cat{A})}-}[rr] & & \lcomod{\cat{L}(\omega)}
\ar@{->}^-{(-){\can{\gG(\cat{A})}}}[rr]
& & \lcomod{\coring{C}}. }
\end{equation*}

Indeed this is a conceptual framework that allows us to compare
certain categories of $\Bbbk$-linear functors with the categories
of comodules over some corings (or left bialgebroids). For
instance, take $R=\Bbbk$ to be a field and $A$ a finite
dimensional $\Bbbk$-algebra. Consider the associated cochain complex
$Q_{\bullet}$ mentioned above and the monoidal
$\Bbbk$-linear category $\Bbbk(\mathbb{N})$ generated by the
natural number $\mathbb{N}$. There is a fiber functor $\chi:
\Bbbk(\mathbb{N}) \to \rmod{\Bbbk}$ defined by $\chi(n)=Q_n$ on
objects and sending the morphism $n \mapsto n+1$ to the
differential $\partial: Q_n \to Q_{n+1}$, for every
$n\in\mathbb{N}$. Using the previous arguments and notations, we
then arrive to the following composition
\begin{equation}\label{Eq:composition}
\xymatrix@C=80pt{ \chain{\Bbbk} \ar@{->}_{\cong}^-{\oO}[r] & \lmod{\gG(\Bbbk(\mathbb{N}))} \ar@{->}^-{Q\tensor{\gG(\Bbbk(\mathbb{N}))}-}[r]
&  \lcomod{\cat{L}(\chi)}
\ar@{->}^-{(-)_{\can{\gG(\Bbbk(\mathbb{N}))}}}[r] &
\lcomod{\lL(A)} }
\end{equation}
where $\oO$ is the canonical equivalence between chain complexes of
$\Bbbk$-vector spaces and left unital $\gG(\Bbbk(\mathbb{N}))$-modules. This in fact is exactly the functor used by D.
Tambara in the proof of \cite[Theorem 4.4]{Tambara}. However,  the above process of  constructing this functor, is actually entirely different from the one presented in \cite{Tambara}. The detailed construction of the functors involved in \eqref{Eq:composition},  as well as  conditions on the extension $ R \to A$ under which this composition gives a monoidal
equivalence form a part of the main aim of this paper.

\subsection{A brief description of the main results.}
Let $\Bbbk$ be a commutative base ring with $1$. Fix a morphism of
$\Bbbk$-algebras $R \to A$. Assume that ${}_RA$ is finitely
generated and projective left $R$-module with a finite dual basis
$\{e_i, {}^*e_i\}_i$. We consider the monoidal functor $-\times_R
A: \bimod{\Re} \to \bimod{R}$, where $(-\times_R-)$ is the
Sweedler-Takeuchi's product \cite{Sweedler:1975b, Takeuchi:1977b}, see the forthcoming subsection. We  obtain that the restriction of
this functor to the category of $\Re$-rings (i.e. the category of
monoids in $\bimod{\Re}$) admits a left adjoint which we denoted
by $\lL: \ring{R} \to \ring{\Re}$. Then $\lL(A)$ the image
of $A$ by the functor $\lL$,   admits a structure of left $R$-bialgebroid
(termed a \emph{coendomorphism bialgebroid}) such that $A$ is a
left $\lL(A)$-comodule ring, see Proposition
\ref{pro: bialgebra} and Corollary \ref{A-comod}. Explicitly, the underlying $\Bbbk$-module $\lL(A)$ is given by the
following quotient of the tensor $\Re$-ring of
$A\tensor{\Bbbk}{}^*A$:
\begin{equation}\label{0-Lm}
\lL(A): = \frac{\tT_{\Re}\lr{A\tensor{}{}^*A}
}{\left\langle\underset{}{}
\sum_i(a\tensor{}e_i\varphi)\tensor{\Re}(a'\tensor{}{}^*e_i) -
(aa'\tensor{}\varphi),\, (1\tensor{}\varphi) -
1\tensor{}\varphi(1)^o \right\rangle_{ \{a,\, a'\, \in A,\,\,
\varphi \,\in {}^*A\} }}.
\end{equation}

On the other hand, we consider the augmented cochain complex of the universal
differential graded ring:
\begin{equation}\label{0-Q}
\xymatrix{Q_{\bullet}: R \ar@{->}^-{1}[r] & A \ar@{->}^-{\partial}[r] & K \ar@{->}^-{\partial_2}[r] & K\tensor{A}K \ar@{->}^-{\partial_3}[r] &
K\tensor{A}K\tensor{A}K \ar@{->}[r] & \cdots \cdots   }
\end{equation}
where $K$ denotes the kernel of $A\tensor{R}A \to A$, the
multiplication of $A$. We check that this is in fact a cochain
complex of left $\lL(A)$-comodules whose components are finitely
generated and projective left $R$-modules. This leads to a fiber
functor $\chi: \Bbbk(\mathbb{N}) \to \bimod{R}$ defined in the
obvious way, as well as to a canonical map $\can{B}:
Q\tensor{B}{}^{\vee}Q \to \lL(A)$, where $B=\Bbbk^{(\mathbb{N})}
\oplus \Bbbk^{(\mathbb{N})}$ (direct sum of copies of $\Bbbk$) is the  ring with enough orthogonal
idempotents attached to the small category $\Bbbk(\mathbb{N})$ (this is the $\Bbbk$-linear category generated by the set of natural numbers $\mathbb{N}$).
Using the fiber functor $\chi$, we first endow   $Q\tensor{B}{}^{\vee}Q$ with a structure of left $R$-bialgebroid, and then
show that $\can{B}$ is an isomorphism of left $R$-bialgebroids. This means, in the sense of \cite{El Kaoutit/Gomez:2004}, that $Q$ is actually a Galois object in the category of left comodules.  In
this way we arrive to our first main result stated below as
Theorem \ref{equivalencia-B}:

\renewcommand{\thetheorem}{{\bf A}}
\usecounter{theorem}
\begin{theorem}
Let $R \to A$ be a $\Bbbk$-algebra map with $A$ finitely generated
and projective as left  $R$-module. Consider the associated left
$R$-bialgebroid $\lL(A)$ (see equation \eqref{0-Lm} above) and
the cochain complex $Q_{\bullet}$ of equation \eqref{0-Q} with its
canonical right unital $B$-action and left $\lL(A)$-coaction,
where $B=\Bbbk^{(\mathbb{N})} \oplus \Bbbk^{(\mathbb{N})}$. Then
the following statements are equivalent
\begin{enumerate}
 \item The right $R$-module ${}_{1\tensor{\Bbbk}R^o}\lL(A)$ is  flat and  the functor $Q\tensor{B}-: \lmod{B} \longrightarrow \lcomod{\lL(A)}$
is an equivalence of monoidal categories;
\item $Q_B$ is a faithfully flat unital module.
\end{enumerate}
\end{theorem}
Obviously, the category of left unital $B$-module $\lmod{B}$ is
isomorphic to the category of  chain complexes of
$\Bbbk$-modules. Thus, Theorem \textbf{A} allows one  freely to transfer  the monoidal model structure  described in  \cite[\S 2.3]{Hovey:Model Catg} for chain complexes of $\Bbbk$-modules,  to the categories of left comodules over coendomorphism bialgebroids. This suggests that certain categories of comodules  can be equipped
with a (monoidal) model structure. Indeed this is one of the main motivations to further investigate categories of comodules over this class of bialgebroids.

Clearly the unit map $\Bbbk \to R$ can be extended to a morphism
of rings with the same set of orthogonal idempotents:
$B=\Bbbk^{(\mathbb{N})} \oplus \Bbbk^{(\mathbb{N})} \to
R^{(\mathbb{N})} \oplus R^{(\mathbb{N})}=C$.  By \cite{El Kaoutit:2006}, this enables us to
consider the usual adjunction between the scalars-restriction
functor and the tensor product functor and, in particular, to
define a canonical map $\can{C}$ with codomain a suitable quotient
of $\lL(A)$. Thus one can try to extend Theorem \textbf{A} to
left unital $C$-modules. In this way we arrive to our second main
theorem which is stated below as Theorem \ref{equivalencia-C}:

\renewcommand{\thetheorem}{{\bf B}}
\usecounter{theorem}
\begin{theorem}
Let $R \to A$ be a $\Bbbk$-algebra map with $A$ finitely generated
and projective as left  $R$-module. Consider $\lL(A)$ the
associated left $R$-bialgebroid (see equation \eqref{0-Lm} above)
 and $\jJ$ the coideal of $\lL(A)$ generated by the set of elements $\{1_{\lL(A)}(r\tensor{}1^o - 1\tensor{}r^o)\}_{r \in R}$;
denote by $\bara{\lL(A)}=\lL(A)/\jJ$ the corresponding
quotient $R$-coring. Consider the cochain complex $Q_{\bullet}$ of
equation \eqref{0-Q} with its structures of right unital
$C$-module and left $\bara{\lL(A)}$-comodule. Then the following
statements are equivalent
\begin{enumerate}
 \item The right $R$-module ${}_{1\tensor{\Bbbk}R^o}\bara{\lL(A)}$ is  flat and  the functor
$Q\tensor{C}-: \lmod{C} \longrightarrow \lcomod{\bara{\lL(A)}}$
is an equivalence of categories; \item $Q_C$ is a faithfully flat unital
module.
\end{enumerate}
\end{theorem}

The problem of obtaining an equivalence of categories as above, is
then closely linked to the faithfully flat  condition on the right
unital module $Q$. This  is in fact not at all easy to check. Our
third main result, which is a combination of Theorem
\ref{Q_B-flat} and Proposition \ref{Q-ff}, gives certain homological
conditions under which $Q$ becomes  flat (or faithfully flat).

\renewcommand{\thetheorem}{{\bf C}}
\usecounter{theorem}
\begin{theorem}
The notations and assumptions are that of Theorem {\bf B}. Assume
further that $A_R$ is finitely generated and projective, and  the
cochain complex $Q_{\bullet}$ is exact and splits, in the sense
that, for every $m \geq 1$, $Q_m\,=\, \partial Q_{m-1} \oplus
\bara{Q}_{m}\,=\,{\rm Ker}(\partial) \oplus \bara{Q}_m$ as right
$R$-modules, for some right $R$-module $\bara{Q}_m$.  Then $Q_C$
is a flat module. In particular, $Q_C$ is faithfully flat in
either one of the following cases.
\begin{enumerate}
 \item $A=R\oplus R t$ ($t^2=0$), the trivial extension of $R$.
 \item $\Bbbk$ is a field and $R$ is a division $\Bbbk$-algebra.
\end{enumerate}
\end{theorem}

As a consequence of Theorems {\bf B} and {\bf C}, we get that for
every $\Bbbk$-algebra $R$, there is a left $R$-bialgebroid $\lL$
such that the category of  chain complexes of left
$R$-modules is equivalent to the category of left comodules over
an epimorphic image of $\lL$. In particular, if $R$ is
commutative, then this equivalence is in fact a monoidal
equivalence.

\subsection{Basic notions and notations.}
Given any Hom-set category $\cat{C}$, the notation $X \in \cat{C}$
means that $X$ is an object of $\cat{C}$. The identity morphism of
$X$ will be denoted by $X$ itself. The set of all morphisms $f:X
\to X'$ in $\cat{C}$ is denoted by ${\rm Hom}_{\cat{C}}(X,\,X')$.
The identity functor of $\cat{C}$ is denoted by $id_{\cat{C}}$. We
denote the dual (or opposite) category of $\cat{C}$ by
$\cat{C}^o$. The class of all natural transformations between two
functors $F$ and $G$ is denoted by ${\rm Nat}(F,G)$.

We work over a ground commutative ring with $1$ denoted by
$\Bbbk$. Up to Section \ref{Main Results}, all rings under
consideration are $\Bbbk$-algebras, and morphisms of rings are
morphisms of $\Bbbk$-algebras. Modules are assumed to be unital
modules and bimodules are assumed to be central $\Bbbk$-bimodules.
For every ring $R$, these categories are denoted by $\lmod{R}$
(left modules), $\rmod{R}$ (right modules) and $\bimod{R}$
(bimodules) respectively. The tensor product over $R$, is denoted
as usual by $-\tensor{R}-$. The unadorned symbol $\otimes$ stands for $\otimes_{\Bbbk}$ the tensor product over $\Bbbk$.

We denote by $\mathit{Ch}(R)$ the category of \emph{chain complexes} of
left $R$-modules. That is, complexes of left modules of the form:
$$
\xymatrix@C=20pt{ (M_{\bullet}, d_{\bullet}): \,\, \cdots \ar@{->}[r] & M_n
\ar@{->}^{d_n}[r] & \cdots \ar@{->}[r] & M_2 \ar@{->}^{d_2}[r] &
M_1 \ar@{->}^{d_1}[r] & M_0\ar@{->}^{d_0}[r] & M_{-1}\ar@{->}[r] &
\cdots \ar@{->}[r] & M_{-n}\ar@{->}^{d_{-n}}[r]&\cdots. }
$$
Let $\chain{R}$ denote the full subcategory of
$\mathit{Ch}(R)$ consisting of positive chain complexes i.e. of
complexes of the form: $$ \xymatrix@C=20pt{ (M_{\bullet}, d_{\bullet}): \,\,
\cdots \ar@{->}[r] & M_n \ar@{->}^{d_n}[r] & \cdots \ar@{->}[r] &
M_2 \ar@{->}^{d_2}[r] & M_1 \ar@{->}^{d_1}[r] & M_0. }
$$

From now on, chain complex of left $R$-modules will stands for an object of the category
$\chain{R}$.
When $R$ is commutative (i.e. commutative $\Bbbk$-algbera), we will considered this category  in a standar way
as a monoidal category with unit object the chain complex $R[0]_{\bullet}$, where $R[0]_0=R$, and $R[0]_n=0$, for $n>0$.

Given an $R$-bimodule $X$, its
$\Bbbk$-submodule of $R$-invariant elements is denoted by
$$ X^R\,\,:=\,\, \left\{\underset{}{} x \in X|\, xr\,=\,rx, \, \forall \, r \in R \right\}. $$ This in fact defines a functor
$ (-)^R: \bimod{R} \to \rmod{\cat{Z}(R)}$, where $\cat{Z}(R)$ is
the centre of $R$. As usual, we use the symbols ${\rm
Hom}_{R-}(-,-)$, ${\rm Hom}_{-R}(-,-)$ and ${\rm Hom}_{R-R}(-,-)$
to denote the Hom-functor of left $R$-linear maps, right
$R$-linear maps and $R$-bilinear maps, respectively.

For two bimodules ${}_RP_S$ and ${}_RQ_S$ over rings $R$ and $S$, we will consider the $\Bbbk$-modules
of $R$-linear maps ${\rm Hom}_{R-}(P,Q)$ as an $S$-bimodule with actions:
$$
s f: p \mapsto f(ps),\, \text{ and }\, f s':
p\mapsto f(p)s',\,\text{ for every } f \in {\rm Hom}_{R-}(P,Q),\, s,
s' \in S,\,\text{ and } p \in P.
$$
Similarly, ${\rm Hom}_{-S}(P,Q)$ is considered as an $R$-bimodule with actions:
$$
r g: p \mapsto rg(p),\, \text{ and }\, g r':
p\mapsto g(r'p),\,\text{ for every } g \in {\rm Hom}_{-S}(P,Q),\, r,
r' \in R,\,\text{ and } p \in P.
$$
Under these considerations, the left dual ${}^*X={\rm Hom}_{R-}(X,R)$  of an $R$-bimodule $X$, is
an $R$-bimodule, as well as its right dual $X^*={\rm Hom}_{-R}(X,\, R)$.

Let $R$ be a ring,  for any $r \in R$, we denote by $\op{r}$ the
same element regarded as an element in the opposite ring $\op{R}$.
Let $R^{\Sf{e}}:=\env{R}$ be the enveloping ring of $R$. Next, we  recall the Sweedler-Takeuchi's  \cite{Sweedler:1975b, Takeuchi:1977b} product on the category of $\Re$-bimodules, usually denoted by $-\times_R-$.  So, given an
$\Re$-bimodule $M$, the underlying $\Bbbk$-module $M$ admits
several structures of $R$-bimodule. Among them, we will select in the forthcoming step  the
following two ones. The first structure is that of the opposite
bimodule ${}_{1\tensor{}R^o}M_{1\tensor{}R^o}$ which we denote by
$\op{M}$. That is, the $R$-biaction on $\op{M}$ is given by
\begin{equation}\label{M0}
r\, \op{m} \,\,=\,\, \op{\lr{m\, (1\tensor{}\op{r})}},\quad \op{m}\,s \,\,=\,\, \op{\lr{(1\tensor{}\op{s})\,m}},\qquad \op{m} \in \op{M},\, r,s \in R.
\end{equation}
The second structure is defined by the left $\Re$-module ${}_{\Re}M$.
That is, the $R$-bimodule $M^{l}\,=\, {}_{R\tensor{}\op{1}}M_R$ whose $R$-biaction is defined by
\begin{equation}\label{LReM}
r \, m^l \,\,=\,\, \lr{(r\tensor{}\op{1}) m}^l,\quad m^l \, s \,\,=\,\, \lr{(1\tensor{}\op{s}) m}^l,\quad m^l \in M^l,\, r,s \in R.
\end{equation}
Now, given $M$ and $N$ two $\Re$-bimodules, we set
$$ \cross{M}{N} \,\,:=\,\,  \lr{{\bf M}\,\tensor{R}\, {\bf N}}^R ,
$$ where ${}_{R}{\bf M}_{R}=M^o$ and ${}_{R}{\bf N}_{R}={}_{R\tensor{}1^o}N_{R\tensor{}1^o}$. The  elements of $M\times _R N$ are   denoted by $\sum_i m_i\times_Rn_i$, for $m_i \in M$ and $n_i \in N$. Henceforth, using these notations and given an element  $\cross{m}{n} \in \cross{M}{N}$, we have the following equalities
\begin{equation}\label{cross-1}
 \cross{\lr{m\, (1\tensor{}r^o)}}{n} \,\,=\,\, \cross{m}{n \, (r\tensor{}1^o)},\text{ and }
\cross{\lr{(1\tensor{}r^o) \,m}}{n} \,\,=\,\, \cross{m}{(r\tensor{}1^o) \, n},
\end{equation}
for every $r, s \in R$. The $\Bbbk$-module $ \cross{M}{N}$ is actually an $\Re$-bimodule with actions:
\begin{equation}\label{Re-MN-1}
 (p\tensor{}q^o) \,\, .  (\cross{m}{n}).\, \, (r\tensor{}s^o) \,\,:=\,\, \cross{\lr{(p\tensor{}1^o) \, m\,
(r\tensor{}1^o) }}{\lr{(1\tensor{}q^o)\,n\, (1\tensor{}s^o))}},
\end{equation}
for every $r, s, p, q \in R$ and $\cross{m}{n} \in \cross{M}{N}$.

On the other hand, since we have $M^o_R=M^l_R$ for every $\Re$-bimodule $M$,
there is a canonical natural transformation (injective at least as $\Bbbk$-linear map)
\begin{equation}\label{Theta}
\xymatrix{\Theta_{M,\,N}: \cross{M}{N} \ar@{->}[rr] & & M^l\tensor{R}N^l.}
\end{equation}

With this product, the $\Re$-bimodule $\cross{S}{T}$ is an $\Re$-ring whenever $S$ and $T$ are. The multiplication of $\cross{S}{T}$ is  defined
componentwise, and the identity element is given by $\cross{1_S}{1_T}$.

An \emph{$R$-ring} $S$ is a monoid in the monoidal category of $R$-bimodules, equivalently,  a $\Bbbk$-algebra map $R \to S$. Dually, an $R$-\emph{coring} is a
comonoid in $\bimod{R}$, which is by definition  a three-tuple
$(\fk{C},\Delta,\varepsilon)$ consisting of $R$-bimodule $\fk{C}$
and two $R$-bilinear maps $\Delta: \fk{C} \to
\fk{C}\tensor{R}\fk{C}$ (comultiplication), $\varepsilon: \fk{C}
\to R$ (counit) satisfying the usual coassociativity and counitary
constraints. In contrast with coalgebras, corings admit several
convolution rings. For instance, the \emph{right convolution} of
an $R$-coring $\coring{C}$, is the right dual $R$-bimodule
$\coring{C}^*$ whose multiplication is defined by
$$ \sigma \,.\, \sigma' \,\, =\,\, \sigma \circ (\sigma'\tensor{R}\coring{C}) \circ \Delta,$$ for all $\sigma, \sigma' \in \coring{C}^*$, and its
unit is the counit $\varepsilon$ of $\coring{C}$. A \emph{morphism
of $R$-corings} is an $R$-bilinear map $\phi: \coring{C} \to
\coring{C}'$ such that $\Delta' \circ \phi = (\phi\tensor{R}\phi)
\circ \Delta$ and $\varepsilon' \circ \phi = \varepsilon$. A left
$\fk{C}$-comodule is pair $(N, \lambda_N)$ consisting of left
$R$-module $N$ and left $R$-linear map $\lambda_N: N \to
\fk{C}\tensor{R}N$ (coaction) compatible in the canonical way with
comultiplication and counit. A morphism of left $\fk{C}$-comodules
is a left $R$-linear map which is compatible with coactions. We
denote by $\lcomod{\fk{C}}$ the category of left
$\fk{C}$-comodules. Right comodules are similarly defined. Given
any morphism of $R$-corings $\phi: \coring{C} \to \coring{C}'$ one
can define, in the obvious way,  a functor $(-)_{\phi}:
\lcomod{\coring{C}} \to \lcomod{\coring{C}'}$ refereed to as the
\emph{coinduction functor}.

For more information on comodules as well as the definitions of
bicomodules and cotensor product over corings, the reader is
referred to \cite{Brzezinski/Wisbauer:2003}. For the notions of
bialgebroids and their basic properties, the reader is referred to
\cite{Bohm:arX08}.

We will also  consider here  rings  with
enough orthogonal idempotents. These  are central $\Bbbk$-modules
$B$ with internal multiplication which  admit a decomposition of
$\Bbbk$-modules $B=\oplus_{p \, \in \cat{P}}B1_p = \oplus_{p \,
\in \cat{P}}1_pB$, where $\{1_p\}_{p \,\in \cat{P}} \subsetneq B$
is a set of orthogonal idempotents. Module over a ring  with
enough orthogonal idempotents stands for $\Bbbk$-central and
\emph{unital module}. Recall that $M$ is a left unital $B$-module
provided that $M$ has an associative left $B$-action  which
satisfies $M=\oplus_{p \, \in \cat{P}} 1_p M$. We denote by
$\lmod{B}$ the category of left unital $B$-modules.

\renewcommand{\thetheorem}{\arabic{section}.\arabic{subsection}.\arabic{theorem}}
\usecounter{theorem}

\section{Coendomorphism and comatrices bialgebroids.}\label{sec:1}

\smallskip

\subsection{Coendomorphism bialgebroid and $\times_R$-comodules.}\label{subsec:1}
In this subsection we  recall the construction of  coendomorphism bialgebroids attached to any finitely generated and projective extension of rings.  We also recall from \cite{Schauenburg:1998} the monoidal structure of the category of comodules over the underlying coring of a given left bialgebroid.
\smallskip


A $\times_R$-coalgebra is an $\Re$-bimodule $\Sf{C}$ together with two $\Re$-bilinear maps $\Delta: \Sf{C} \to \cross{\Sf{C}}{\Sf{C}}$ (comultiplication) and  $\varepsilon : \Sf{C} \to \End{\Bbbk}{R}$ (counit) which satisfy the coassociativity and counitary properties in the sense of  \cite[\S 4, Definition 4.5]{Takeuchi:1977b}, see also \cite{Brzezinski/Militaru:2002} and \cite{Schauenburg:1998}.  A $\times_R$-coalgebra $\Sf{C}$ is said to be an $\times_R$-\emph{bialgebra} provided that comultiplication and
counit are morphisms of $\Re$-rings.

A left $\times_R$-$\Sf{C}$-\emph{comodule}, is a pair
$(X,\lambda_{X})$ consisting of an $R$-bimodule $X$ and an
$R$-bilinear map $\lambda_X: X \to \cross{\Sf{C}}{X}$ satisfying,
in some sense, the coassociativity and
counitary axioms. Morphism between left
$\times_R$-$\Sf{C}$-comodules are $R$-bilinear maps compatible in
the obvious way with the left $\times_R$-$\Sf{C}$-coactions. This
leads to the definition of the category of left
$\times_R$-$\Sf{C}$-comodules. When $\Sf{C}$ is a
$\times_R$-bialgebra, this category becomes a monoidal category
\cite[Proposition 5.6]{Schauenburg:1998}, and the forgetful
functor to the category of $R$-bimodules is a monoidal functor.
There is a strong relation which will be clarified in the sequel,  between the category of left
$\times_R$-comodules over an $\times_R$-bialgebra and the category
of left comodules over the  underlying $R$-coring whose structure
maps are
$$
\xymatrix@C=40pt{ \Sf{C} \ar@{->}^-{}[r] & \cross{\Sf{C}}{\Sf{C}} \ar@{->}^-{\Theta_{C,C}}[r] & \Sf{C}^l\tensor{R}\Sf{C}^l,} \quad
\xymatrix@C=60pt{\Sf{C} \ar@{->}^-{\varepsilon(-)(1_R)}[r] & R,}
$$ where $\Theta_{-,-}$ is the  natural transformation of equation \eqref{Theta}.

Let $A$ be an $R$-ring, that is,  a $\Bbbk$-algebra map $R \to A$, and denote by ${}^*A$ the dual of the $\Bbbk$-module $A_{\Bbbk}$, i.e. ${}^*A\,=\, \hom{\Bbbk}{A}{\Bbbk}$ .
We consider the tensor product $A\tensor{}{}^*A$ as an $\Re$-bimodule in the following way
\begin{equation}\label{Eq:DagaM}
(r\tensor{}s^o) \, .\, (a\tensor{}\varphi) \, .\, (p\tensor{}q^o)\,\,=\,\, (rap)\tensor{}(q\varphi s), \quad p,q,r,s \, \in R, \text{ and } a, \in A, \varphi \in {}^*A,
\end{equation}
where $A$ and ${}^*A$ are considered as $R$-bimodules in the usual way.

Assume that ${}_RA$ is finitely generated and projective module, and fix a left dual basis
$\{e_i,{}^*e_i\}_{1 \leq i \leq n}$. Define
the $\Re$-ring $\lL(A)$  by the quotient algebra
\begin{equation}\label{lmC}
\lL(A)\,=\, \tT_{\Re}\lr{A\tensor{}{^*A}} / \iI
\end{equation}
where $\tT_{\Re}\lr{A\tensor{}{}^*A}\,=\, \underset{n
\in\mathbb{N}}{\bigoplus} \lr{A\tensor{}{}^*A}^{\underset{\Re}{\otimes} {} ^n}$ is the tensor $\Re$-ring of the $\Re$-bimodule $A\tensor{}{^*A}$ and where  $\iI$ is the two-sided  ideal generated by the set
\begin{equation}\label{iI_C}
\left\{ \underset{}{}
\sum_i\lr{(a\tensor{}e_i\varphi)\tensor{\Re}(a'\tensor{}
{}^*e_i)} \,-\, (aa'\tensor{}\varphi);\,\,
1_R\tensor{}\varphi(1_A)^o\,-\, (1_A\tensor{}\varphi)
\right\}_{a,\,a' \in\, A,\,\, \varphi \in \,{ }^*A}.
\end{equation}

We denote by
\begin{equation}\label{Eq:pi}
\pi_A: \tT_{\Re}(A\tensor{}{}^*A) \to \lL(A)
\end{equation}
the canonical
projection. From now on, given a homogeneous elements
$(a\tensor{}\varphi) \in \tT_{\Re}(A\tensor{}{}^*A)$
of degree one, we
denote by $\pi_A(a\tensor{}\varphi)$ its image in $\lL(A)$.
Next, we recall the structure of $\times_R$-bialgebra of the object $\lL(A)$, which is denoted by $a_R(A,A)$ in \cite[Remark 1.7]{Tambara}. The underlying  structure of an $\Re$-ring, is given by the following
composition of algebra maps
$$ \xymatrix{ \Re \ar@{->}^-{\iota_0}[r] & \tT_{\Re}(\lL(A)) \ar@{->}^-{\pi_A}[r] & \lL(A), }$$ where $\iota_{m}$ denotes
the canonical $\Re$-bilinear injection in degree $m\geq 0$.

\begin{proposition}\label{pro: bialgebra}
Let $A$ be an $R$-ring which is finitely generated and projective
as left $R$-module with dual basis $\{({}^*e_i,e_i)\}_i$. Then
$\lL(A)$ is a $\times_R$-bialgebra with structure maps
$$\xymatrix@R=0pt@C=50pt{\lL(A) \ar@{->}^-{\Delta}[r] & \lL(A) \times_R \lL(A),
\\ \pi_A(a\tensor{}\varphi) \ar@{|->}[r] & \sum_j \pi_A(a\tensor{}{}^*e_j)\times_R \pi_A(e_j\tensor{}\varphi) }
\quad \xymatrix@R=0pt@C=40pt{ \lL(A) \ar@{->}^-{\varepsilon}[r] & \mathrm{End}_{\Bbbk}(R) \\ \pi_A(a\tensor{}\varphi) \ar@{|->}[r]
& \left[\underset{}{} r \mapsto \varphi(ar) \right].  }$$
\end{proposition}

The  relation between the $R$-ring structure of $A$
and the $\times_R$-bialgebra structure  of $\lL(A)$, is expressed as follows.

\begin{corollary}\label{A-comod}
Let $A$ be an $R$-ring such that ${}_RA$ is finitely generated and
projective and $\lL(A)$ the associated $\cross{}{}$-bialgebra
defined in Proposition \ref{pro: bialgebra}. Then $A$ is a left
$\times_R$-$\lL(A)$-comodule $R$-ring, that is, $A$ admits a
left $\times_R$-$\lL(A)$-coaction
$$\lambda_A: A \longrightarrow \cross{\lL(A)}{A},\quad \lr{a \longmapsto \sum_j\cross{\pi_A(a\tensor{}{}^*e_j)}{e_j}} $$
which is also a morphism of $R$-rings.
\end{corollary}

The $\times_R$-bialgebra $\lL(A)$ defined in Proposition
\ref{pro: bialgebra} is refereed to as \emph{coendomorphism
$R$-bialgebroid} since by \cite[Theorem
3.1]{Brzezinski/Militaru:2002}, $\lL(A)$ is in fact a (left)
bialgebroid whose structure of $\Re$-ring is the map
$$\pi_A \circ \iota_0: \Re \longrightarrow \lL(A),$$ and its structure of $R$-coring is given as follows. The underlying $R$-bimodule is
$\lL(A)^l={}_{\Re}\lL(A)$, the comultiplication and counit are given by
\begin{equation}\label{Delta-l}
\Delta: \lL(A)^l \longrightarrow  \lL(A)^l \tensor{R}\lL(A)^l, \quad \lr{\pi_A(a\tensor{}\varphi)
\longmapsto \sum_i \pi_A(a\tensor{}{}^*e_i) \tensor{R} \pi_A(e_i\tensor{}\varphi)},
\end{equation}
\begin{equation}\label{epsilon-l}
\varepsilon: \lL(A)^l \longrightarrow R,\quad \lr{ \pi_A(a\tensor{}\varphi) \longmapsto \varphi(a)}.
\end{equation}

\begin{example}\label{main-example}
Let $A=R\oplus R t$  be the trivial generalized $R$-ring, i.e. the
$R$-ring which is free as left $R$-module with basis $1=(1,0)$ and
$\fk{t}=(0,t)$ such that $\fk{t}^2=0$. Using \eqref{lmC} and
Proposition \ref{pro: bialgebra}, we can easily check that
$\lL(A)$ is an $R$-bialgebroid generated by the image of $\Re$
and two $\Re$-invariant elements $\{x,y\}$ subject to the
relations $xy + yx \,=\,0$, $x^2 \,=\, 0$. The comultiplication
and counit of it underlying $R$-coring are given by
\begin{eqnarray*}
  \Delta(x) &=& x\tensor{R}1 + y\tensor{R}x,\quad \varepsilon(x)\,\,=\,\,0 \\
  \Delta(y) &=& y\tensor{R}y,\quad \varepsilon(y)\,\,=\,\,1.
\end{eqnarray*}
The ring $A$ is a  left $\lL(A)$-comodule ring with coaction: $\lambda: A\to \lL(A)\tensor{R}A$ sending
$$\lambda(1_A)\,=\, 1_{\lL(A)}\tensor{R}1_A,\quad \lambda(\fk{t})\,=\, x\tensor{R}1_A + y\tensor{R}\fk{t},$$  extended by $R$-linearity to the whole set of elements of $A$.
\end{example}


In \cite{Schauenburg:1998} it was
shown that the category of left $\times_R$-comodules over an
$\times_R$-bialgebra is a monoidal category such that the
forgetful functor to the category of $R$-bimodules is a monoidal
functor. What we will need in the sequel is a monoidal structure
on the category of left $\lL(A)$-comodules where $\lL(A)$ is viewed as
an $R$-coring with structure maps \eqref{Delta-l} and
\eqref{epsilon-l}.  The following lemma is a consequence of
\cite[Proposition 5.6]{Schauenburg:1998}, see also \cite[3.6]{Bohm:arX08}.
\begin{lemma}\label{Mono-comod-1}
Let $\lL$ be any left $R$-bialgebroid. Then the category of left
$\times_R$-$\lL$-comodule is isomorphic to the category of left
$\lL^l$-comodules over the underlying $R$-coring $\lL^l$. In
particular, the category of left $\lL^l$-comodules inherits a
monoidal structure with unit object $(R, R\to \lL^l)$ and the left
forgetful functor $U: \lcomod{\lL^l} \to \lmod{R}$ factors
throughout a  monoidal functor into the category of $R$-bimodules.
Thus, we have a commutative diagram
$$
\xymatrix{ \lcomod{\lL^l} \ar@{->}^-{U}[rr] \ar@{-->}[rd] & & \lmod{R} \\ & \bimod{R}\ar@{->}[ru] }
$$
where the dashed arrow is a monoidal functor.
\end{lemma}

Summing up, given two left $\lL^l$-comodules $(X,\lambda_X)$ and $(Y,\lambda_Y)$,
using Lemma \ref{Mono-comod-1}, we can consider $(X\tensor{R}Y,\lambda_{X\tensor{R}Y})$ as a
left $\lL^l$-comodule with coaction
\begin{equation}\label{multiplication-Comod}
\lambda_{X\tensor{R}Y}: X\tensor{R}Y \to \lL^l\tensor{R} X\tensor{R}Y, \quad \lr{x\tensor{R}y \longmapsto \sum_{(x),(y)} (x_{(-1)}y_{(-1)})^l
\tensor{R} (x_{(0)}\tensor{R}y_{(0)}), }
\end{equation}
where we have considered $X$ as $R$-bimodule with the right
$R$-action given by the action
\begin{equation}\label{right}
 x\,r\,\,=\,\, \sum_{(x)} \varepsilon\lr{x_{(-1)} (r\tensor{}1^o)} x_{(0)},\, \text{ for every  }\, r \in R \text{ and } x  \in X.
\end{equation}

\subsection{The complex of left $\lL$-comodules $Q_{\bullet}$. }\label{cochain-Q}
Keep the assumptions and notations of 
subsection \ref{subsec:1}, that is,
we are considering an $R$-ring $A$ over a fixed $\Bbbk$-algebra $R$.
Let us denote by
$$ K = \mathrm{Ker}\lr{\xymatrix{A\tensor{R}A \ar@{->}^-{\mu}[r] & A}}$$
the kernel of the multiplication $\mu$ of $A$ with canonical
derivation
$$\xymatrix@R=0pt{ A \ar@{->}^-{\partial}[r] & K \\ a \ar@{|->}[r] & \lr{\partial a\,\,=\,\,1\tensor{R}a -a \tensor{R}1}. }$$
The associated cochain complex is denoted by
$$\xymatrix{Q_{\bullet}: R \ar@{->}^-{\partial_0=1}[r] & A \ar@{->}^-{\partial_1=\partial}[r] & K \ar@{->}^-{\partial_2}[r] & K\tensor{A}K \ar@{->}^-{\partial_3}[r] &
K\tensor{A}K\tensor{A}K \ar@{->}[r] & \cdots \cdots   }$$ where
$\partial_n: Q_n \to Q_{n+1}$ sends $a_0\partial a_1\tensor{A}
\cdots \tensor{A} \partial a_{n-1}$ to $\partial
a_0\tensor{A}\partial a_1\tensor{A}\cdots \tensor{A}\partial
a_{n-1}$, $n \geq 2$.

The following lemma, which will play a key role in Subsection
\ref{Result-3}, characterizes a split ring extension $R \to A$ (in
$\rmod{R}$) in terms of the cochain complex $Q_{\bullet}$.
\begin{lemma}\label{mu-splits}
Let $A$ be any  $R$-ring. Then the following conditions are equivalent.
\begin{enumerate}[(i)]
 \item The unit $u: R \to A$ is a split monomorphism in
$\rmod{R}$.
\item The cochain complex $Q_{\bullet }$ is exact and
splits, in the sense that, for every $m\geq 1$, there is a right $R$-module $\bara{Q}_{m}$ such that $Q_{m}=\partial
Q_{m-1}\oplus \bara{Q}_{m}= {\rm Ker}(\partial) \oplus \bara{Q}_{m}$, as right $R$-modules.
\end{enumerate}
\end{lemma}
\begin{proof}
 $(ii) \Rightarrow (i) $ It is trivial.\\
$(i) \Rightarrow (ii) $. Let us denote by $u^c: A \to \bara{A}$ the cokernel of $u: R \to A$ in $\bimod{R}$.
Put $\Omega_0:=R$, $\Omega_1:=A$, and $\Omega_n:=A\tensor{R}\bara{A}\tensor{R} \cdots \tensor{R}\bara{A}$, $(n-1)$-fold  $\bara{A}$, for $n \geq 2$.
Consider now the following split exact sequence of right $R$-modules
$$
\xymatrix{ 0 \ar@{->}[r] & \bara{A}^{\tensor{R} \, n} \ar@{->}^-{\gamma_n}[r] & A\tensor{R}\bara{A}^{\tensor{R} \, n}
\ar@{->}[r] & \bara{A}^{\tensor{R} \, n+1} \ar@{->}[r] & 0,}
$$
where $\gamma_n= u \tensor{R} \bara{A}^{\tensor{R} \, n} $, for $n \geq 1$. In view of this, we have a split exact cochain complex of right $R$-modules
$$\xymatrix@C=40pt{\Omega_{\bullet} : \Omega_0 \ar@{->}^-{d_0}[r] & \Omega_1 \ar@{->}^-{d_1}[r] & \Omega_2 \ar@{->}^-{d_2}[r] & \Omega_3 \ar@{->}[r] & \cdots , }$$
with differential $d_0=u$, $d_1= \gamma_1 \circ u^c$, $d_n= \gamma_n \circ (u^c\tensor{R}\bara{A}^{\tensor{R} \, n-1})$, for $n \geq 2$.
Since $\Omega_2$ is the cokernel of the map $A\tensor{R}u$, and the later split by $\mu$ the multiplication of $A$,
we obtain the following split exact sequence of $R$-bimodules
$$\xymatrix@C=40pt{ 0 \ar@{->}[r] & A  \ar@{->}^-{A\tensor{R}u}[r] & A\tensor{R}A \ar@{->}^-{A\tensor{R}u^c}[r] & \Omega_2 \ar@{->}[r] & 0.}
$$
This gives the split exact sequence
$$
\xymatrix{ 0 \ar@{->}[r] & \Omega_2   \ar@{->}^-{}[r] & A\tensor{R}A \ar@{->}^-{\mu}[r] & A \ar@{->}[r] & 0.}
$$
Thus we have an $R$-bilinear isomorphism  $\omega_2: \Omega_2 \to Q_2=K$.
Henceforth, there is an unique $A$-bimodule structure on $\Omega _{2}$ which renders $\omega _{2}$ an $A$-bilinear
isomorphism, namely
\begin{equation*}
a\cdot \left( x\otimes _{R}\overline{y}\right) \cdot b=ax\otimes _{R}
\overline{yb}-axy\otimes _{R}\overline{b}, \quad\text{ for every }a,x,y,b\in A,
\end{equation*}
wherein the notation $u^c(z)=\bara{z}$, for every $z \in A$, have been used.
Define iteratively
$\omega _{n}:\Omega _{n}\rightarrow Q_{n}$, for all $n\geq 3,$ as the composition
$$
\xymatrix@C=50pt{\Omega _{n}=\Omega _{n-1}\otimes _{R}\overline{A}\cong \Omega _{n-1}\otimes
_{A}\left( A\otimes _{R}\overline{A}\right) =\Omega _{n-1}\otimes _{A}\Omega
_{2}   \ar@{->}^-{\omega _{n-1}\otimes _{A}\omega _{2}}[r]
& Q_{n-1}\otimes _{A}K=K^{\otimes _{A}n-1}=Q_{n}.}
$$
By construction, $\omega_{\bullet}: \left( \Omega _{\bullet
},d_{\bullet }\right) \to \left( Q_{\bullet },\partial _{\bullet
}\right)$ is a morphism of complexes of $R$-bimodules. We leave to
the reader to check that $\omega_{\bullet}$ is in fact an
isomorphism of cochain complexes. Now, since $\left( \Omega
_{\bullet },d_{\bullet }\right) $ is split exact in right
$R$-modules, then so is $\left( Q_{\bullet },\partial _{\bullet
}\right) $.
\end{proof}

\begin{remark}\label{Remark-mu-splits}
In the case of finitely generated and projective extension of rings, the left version of
condition $(i)$ in Lemma \ref{mu-splits} implies that ${}_RA$ is
in fact faithfully flat module (see, for example \cite[Chap. I,
Proposition 9, page 51]{Bourbaki:Chap. I-IV}). In this case, one
can easily show that $Q_{\bullet}\tensor{R}A$ is homotopically trivial
which by \cite[Th\'eor\`eme 2.4.1]{Godement:1973} gives condition $(ii)$. In this way, Lemma \ref{mu-splits}
can be seen as a generalization of \cite[Propositions 6.1,
6.2]{Artin:1969}.
\end{remark}

The convolution product on the left dual chain complex of
$Q_{\bullet}$ is given as follows: For every $\varphi \in {}^*Q_n$
and $\psi \in {}^*Q_m$ with $n,m \geq 1$, we have a left
$R$-linear map
\begin{equation}\label{1-star}
\xymatrix@R=0pt{ \varphi \star \psi: Q_{n+m} \ar@{->}[r] & R \\ \quad x\tensor{A}\partial(a)\tensor{A}y \ar@{|->}[r] &
\varphi\lr{x\psi(ay)}\,-\, \varphi\lr{xa\psi(y)}, }
\end{equation}
where $x \in Q_n$ , $y \in Q_m$, and $a \in A$. The convolution
product with zero degree element is just the left and right
$R$-actions of ${}^*Q_n$, for every $n \geq 1$, namely
\begin{equation}\label{0-star}
\xymatrix@R=0pt{ r \star \varphi : Q_{n} \ar@{->}[r] & R \\ \quad x \ar@{|->}[r] &
\varphi(x \, r), } \qquad \xymatrix@R=0pt{ \varphi \star s : Q_{n} \ar@{->}[r] & R \\ \quad x \ar@{|->}[r] &
\varphi(x)\, s, }
\end{equation}
for every elements $r, s \in R$ and $\varphi \in {}^*Q_n$.

\begin{remark}\label{Q-normalization}
The convolution product defined in \eqref{1-star} and
\eqref{0-star} derives from the structure of comonoid of the
cochain complex $Q_{\bullet}$ viewed as an object in the monoidal
category of cochain complexes of $R$-bimodules. Precisely, the
identity map $A\tensor{R}\cdots \tensor{R} A=A^{\tensor{R}\,n}
=A^{\tensor{R}\,p} \tensor{R} A^{\tensor{R}\,q}$, for $p+q=n$,
rereads as a map $Q_{n} \to Q_p\tensor{R}Q_q$ sending
$x\tensor{A}\partial a\tensor{A} y \mapsto x\tensor{R}ay
-xa\tensor{R}y$, for every $x \in Q_p$, $a \in A$ and $y \in Q_q$.
Thus $Q\,=\,\oplus_{n \geq 0} Q_n$ has  a structure of
differential $R$-coring in the sense of \cite[pages 6,
7]{Cartan:1958}. Since each $Q_n$ is finitely generated and
projective left $R$-module (see Lemma \ref{dual-basis} below, of course under the same assumption for the left module ${}_RA$), the
comultiplication of $Q$ is transferred to the  graded left dual
${}^{\vee}Q\,=\, \oplus_{n \geq 0} {}^*Q_n$ which gives a
multiplication defined explicitly by \eqref{1-star} and
\eqref{0-star}.
A comonoidal structure on $Q_{\bullet}$ could also be obtained by
transferring some comonoidal structure of the Amitsur cosimplicial
object of $R$-bimodules induced by $A$ (see \cite{Artin:1969}), using for this
the normalization functor and it structure of comonoidal functor
obtained from Eilenberg-Zilber Theorem, see \cite[Theorem 8.1,
Exercise 4. p. 244]{MacLane:Homology}  (of course in their dual
form). It seems that Tambara's approach \cite{Tambara} runs in
this direction. Anyway this approach uses a slightly variant of
the category of  cosimplicial groups endowed with some monoidal
structure which is not the usual one. Since our methods run in a
different way, we will not make use of the normalization process
here.
\end{remark}

In all what follows, we will fix a (left) finitely generated and projective extension $R \to A$
 with dual basis $\{e_i,{}^*e_i\}_{1 \leq i\leq n}$. We will denote by $\lL:=\lL(A)$ the corresponding left $R$-bialgebroid  coming from Proposition \ref{pro: bialgebra}, and by $\pi$
the projection $\pi_A$ defined in  \eqref{Eq:pi}.

Using this dual basis, one can check that ${}_RQ_2={}_RK$
is finitely generated and projective module whose dual basis is
given by the set $\{(e_i\partial e_j, {}^*e_i \star
{}^*e_j)\}_{i,\, j}$. Moreover, we have
\begin{lemma}\label{dual-basis}
Each  $Q_{n}$, $n \geq 0$, is finitely generated and projective as left $R$-module. Furthermore, if $\{(\omega_{n,\alpha},
{}^*\omega_{n,\alpha})\}_{\alpha}$ is a dual basis  for $Q_n$ with $n\geq 1$, then
$\{(\omega_{n,\alpha} \tensor{A}\partial \omega_{m,\beta}, {}^*\omega_{n,\alpha} \star {}^*\omega_{m,\beta} )\}_{\alpha,\,\beta}$
is a dual basis for $Q_{n+m}$, while $\{(\omega_{n,\alpha}\tensor{A}\omega_{m,\beta},
{}^*\omega_{n,\alpha}\star\partial{}^*\omega_{m,\beta})\}_{\alpha,\,\beta}$
is a dual basis for $Q_{n+m-1}$ when $m \geq 2$.
\end{lemma}
\begin{proof}
Straightforward.
\end{proof}
 The cochain complex $Q_{\bullet}$ is actually a complex of left $\lL$-comodules.
\begin{proposition}\label{Q-comod}
The cochain complex $Q_{\bullet}$ is a complex of left
$\lL$-comodules. For $n=0$, the coaction is given by $(R \to \lL,
r \mapsto \pi(r\tensor{}1^o))$ and, for $  n\geq 1$, by
$\lambda_n: Q_n \to  \lL \tensor{R} Q_n$ sending
\begin{footnotesize}
\begin{multline}
\xymatrix{ a_0\partial a_1 \tensor{A}\cdots\tensor{A}\partial a_{n-1}  \ar@{|->}[rr] & & } \\
\underset{i_0,\, i_1,\, \cdots \, ,\,i_{n-1}}{\sum} \pi(a_0\tensor{}{}^*e_{i_0})\cdots\pi(a_{n-1}\tensor{}{}^*e_{i_{n-1}}) \tensor{R}
\lr{e_{i_0}\partial e_{i_1} \tensor{A}\cdots\tensor{A}\partial e_{i_{n-1}}  }.
\end{multline}
\end{footnotesize}
\end{proposition}
\begin{proof}
The statement is trivial for $n=0$. For $n \geq 1$, the
coassociativity of $\lambda_n$ is deduced using Lemma \ref{dual-basis} which assert that $\{(e_{i_0}\partial e_{i_1} \tensor{A}\cdots\tensor{A}\partial
e_{i_{n-1}} ,\, {}^*e_{i_0} \star \cdots \star {}^*e_{i_{n-1}}
\}_{i_0,\,i_1,\, \cdots,i_{n-1}}$ is  a dual basis for $Q_n$. Here  each ${}^*e_{i_0} \star \cdots \star
{}^*e_{i_{n-1}}$ is the $n$-fold convolution product defined in
\eqref{1-star}. The rest of the proof uses the fact that each coaction $\lambda_n$, $n \geq 1$, satisfies the equality
\begin{footnotesize}
\begin{equation}\label{sin-1}
\lambda_n\lr{\partial b_1 \tensor{A}\cdots\tensor{A}\partial b_{n-1} } \,\, = \,\,
\underset{ i_1,\, \cdots\, ,\,i_{n-1}}{\sum} \pi(b_1\tensor{}{}^*e_{i_1})\cdots\pi(b_{n-1}\tensor{}{}^*e_{i_{n-1}}) \tensor{R}
\lr{\partial e_{i_1} \tensor{A}\cdots\tensor{A}\partial e_{i_{n-1}}  }.
\end{equation}
\end{footnotesize}
\end{proof}

The following lemma will be used in the sequel.
\begin{lemma}\label{Q-comod-1}
Given two elements $u_n=a_0\partial a_1\tensor{A}\cdots\tensor{A}\partial a_{n-1} \in Q_n$ and
$u_m=b_0\partial b_1 \tensor{A}\cdots\tensor{A}\partial b_{m-1}
\in Q_m$ with $n,m\geq 1$. Then
\begin{footnotesize}
\begin{multline*}
 \lambda_{n+m-1}(u_n\tensor{A}u_m) \,=\,
\sum_{i_0,\cdots,\,i_{n-1},\,j_0, \cdots,j_{m-1}} \lr{\pi(a_0\tensor{}{}^*e_{i_0})
\cdots\pi(a_{n-1}\tensor{}{}^*e_{i_{n-1}}) \pi(b_0\tensor{}{}^*e_{j_0})
\cdots\pi(b_{m-1}\tensor{}{}^*e_{j_{m-1}})} \\ \tensor{R}
\lr{e_{i_0}\partial e_{i_1} \tensor{A} \cdots \tensor{A}\partial e_{i_{n-1}}
\tensor{A}e_{j_0}\partial e_{j_1} \tensor{A} \cdots \tensor{A}\partial e_{j_{m-1}}}.
\end{multline*}
\end{footnotesize}
Furthermore, for every $u \in Q_n$, $n \geq 1$ and $v \in Q_m$, $m \geq 1$, we have
$$
\lambda_{n+m-1}(u\tensor{A}v)\,\,=\,\, \sum u_{(-1)}v_{(-1)} \tensor{R} (u_{(0)}\tensor{A}v_{(0)}), $$  and
$$
\lambda_{n+m}(u\tensor{A}\partial v)\,\,=\,\, \sum u_{(-1)}v_{(-1)} \tensor{R} (u_{(0)}\tensor{A}\partial v_{(0)}),
$$
where Sweedler's notation for coactions is used.
\end{lemma}
\begin{proof}
The proof of the first claim is based upon the observation that
the coaction of any $Q_k=K\tensor{A} \cdots \tensor{A}K$
($(k-1)$-times), with $k \geq 2$, is induced from that of
$A\tensor{R}\cdots\tensor{R}A$ ($k$-times). The later is a left
$\lL$-comodule, by Corollary \ref{A-comod} and Lemma
\ref{Mono-comod-1}, using the coactions described in
\eqref{multiplication-Comod}. The last statement is deduced from
the first one by left $R$-linearity.
\end{proof}

\subsection{The infinite comatrix bialgebroid induced by $Q_{\bullet}$.}\label{comatrix}

Let $Q_{\bullet}$ be the cochain complex of $\lL$-comodules
considered in Proposition \ref{Q-comod}. In this subsection we
will construct a left bialgebroid associated to $Q_{\bullet}$ and
a canonical map from this left bialgebroid to $\lL$. First we
recall from \cite{El Kaoutit/Gomez:2003,El Kaoutit/Gomez:2004} the
notion of infinite comatrix coring and the canonical map. A
different approach to this notion can be found in
\cite{Wisbauer:2006}, \cite{Caenepeel/Groot/Vercruysse:2006} and
\cite{Gomez/Vercruysse:2005}. We should mention here that this
object coincides with the one already  constructed  in the context
of Tannaka-Krein duality over fields or commutative rings, see
\cite{Deligne:1990}, \cite{Bruguieres:1994},
\cite{Joyal/Street:1991} and \cite{P.H. Hai:2008}, see also
\cite{McCrudden:2000}. However, the description given in \cite{El
Kaoutit/Gomez:2004} in terms of tensor product over a ring with
enough orthogonal idempotents, seems easier to handle from a
computational point of view.

Let $\cat{A}$ be a small full sub-category of an additive
category. Following \cite[page 346]{Gabriel:1962}, we can associate to $\cat{A}$ the
ring with enough orthogonal idempotents $S=\oplus_{\fk{p},\,
\fk{p}' \in \cat{A}}{\rm Hom}_{\cat{A}^{o}}(\fk{p},\,\fk{p}')$,
where $\cat{A}^{o}$ is the opposite category of $\cat{A}$. The
category of left unital $S$-module is denoted by $\lmod{S}$.

Let us denote by $add({}_RR)$ the full sub-category of $\lmod{R}$
consisting of all finitely generated and projective left
$R$-modules.  Let $\chi: \cat{A} \rightarrow add({}_RR)$ be a
faithful functor, refereed to as \emph{fiber functor}. We denote
by $\fk{p}^{\chi}$ the image of $\fk{p} \in \cat{A}$ under $\chi$
or by $\fk{p}$ itself if no confusion arises.  Consider the left
$R$-module direct sum of the $\fk{p}$'s: $\Sigma= \oplus_{\fk{p}
\in \cat{A}} \fk{p}$ (i.e. $\Sigma= \oplus_{\fk{p} \in
\cat{A}}\fk{p}^{\chi}$) and the right $R$-module direct sum of
their duals: ${}^{\vee}\Sigma= \oplus_{\fk{p} \in \cat{A}}
{}^*\fk{p}$. It is clear that ${}^{\vee}\Sigma$ is a left unital
$S$-module while $\Sigma$ is a right unital $S$-module. In this
way $\Sigma$ becomes an $(R,S)$-bimodule and ${}^{\vee}\Sigma$ an
$(S,R)$-bimodule. Then $\Sigma \tensor{S} {}^{\vee}\Sigma$ is now
an $R$-bimodule whose elements are described as a finite sum of
diagonal ones, i.e. of the form $\iota_{\fk{p}}(u_{\fk{p}})
\tensor{S}\iota_{{}^*\fk{p}}(\varphi_{\fk{p}})$ where
$(u_{\fk{p}},\varphi_{\fk{p}}) \in {\fk{p}}^{\chi} \times
({}^*\fk{p}^{\chi})$ and $\iota_{-}$ are the canonical injections
in ${}^{\vee}\Sigma$ and $\Sigma$.  From now on, we will write
$u_{\fk{p}}\tensor{S}\varphi_{\fk{p}}$ instead of
$\iota_{\fk{p}}(u_{\fk{p}}) \tensor{S}
\iota_{{}^*\fk{p}}(\varphi_{\fk{p}})$ to denote a generic element
of $\Sigma \tensor{S} {}^{\vee}\Sigma$.

This bimodule admits a
structure of an $R$-coring given by the following comultiplication
\begin{equation}\label{Delta}
\xymatrix@R=0pt{ \Delta: \Sigma \tensor{S} {}^{\vee}\Sigma
\ar@{->}[r] & \left( \Sigma \tensor{S} {}^{\vee}\Sigma \right)
\tensor{R} \left(\Sigma \tensor{S} {}^{\vee}\Sigma \right) \\
u_{\fk{p}} \tensor{S} \varphi_{\fk{p}}
\ar@{|->}[r] & \sum_i u_{\fk{p}} \tensor{S}
{}^*u_{\fk{p},\, i} \tensor{R}
u_{\fk{p},\,i} \tensor{S} \varphi_{\fk{p}}, }
\end{equation}
where, for a fixed $\fk{p} \in \cat{A}$, the finite set $\{(u_{\fk{p},
\, i},{}^*u_{\fk{p},\, i}) \}_i \subset \fk{p} \times {}^*\fk{p}$ is a left
dual basis of the left $R$-module $\fk{p}$.
The counit is just the evaluating map.
Note that this comultiplication is independent from the chosen
bases. With this structure $\Sigma \tensor{S}{}^{\vee}\Sigma$ is
refereed to as the \emph{infinite comatrix coring} associated to
the small category $\cat{A}$ and the fiber functor $\chi$. On the other
hand, each of the left $R$-modules $\fk{p}^{\chi}$ is actually a
left $(\Sigma\tensor{S}{}^{\vee}\Sigma)$-comodule with coaction,
using the above notation is given by
\begin{equation}\label{td-lambda}
 \td{\lambda}_{ \fk{p}}: \fk{p} \longrightarrow \Sigma\tensor{S}{}^{\vee}\Sigma \tensor{R} \fk{p}, \quad \lr{u \longmapsto
\sum_i u \tensor{S}
{}^*u_{\fk{p},\, i} \tensor{R}\, u_{\fk{p},\,i}
 }.
\end{equation}

Another description of the infinite comatrices is given in
\cite[Proposition 5.2]{El Kaoutit/Gomez:2004} which establishes an
isomorphism of $R$-bimodules
\begin{equation}\label{R(S)}
\Sigma\tensor{B}{}^{\vee}\Sigma \,\, \cong \,\, \frac{\underset{\fk{p} \in\, \cat{A}}{\oplus} \,\fk{p}\tensor{T_{\fk{p}}}{}^*\fk{p}}{\left\langle
u \fk{t} \tensor{T_{\fk{q}}} \varphi \, - \,
u  \tensor{T_{\fk{p}}} \fk{t} \varphi
\right\rangle}_{\left\{ u \,\in \,\,\fk{p},\,\, \varphi\, \in\,\, {}^*\fk{q}, \,\,\fk{t}\, \in \, \, T_{\fk{q},\fk{p}} \right\} }
\end{equation}
where $T_{\fk{p}}:={\rm End}_{\cat{A}^o}(\fk{p})$  and
$T_{\fk{p},\,\fk{q}}\,:=\, \hom{\cat{A}^o}{\fk{p}}{\fk{q}}$, for
every objects $\fk{p} , \fk{q}$ in $\cat{A}$.

Now, let $\coring{C}$ be an $R$-coring and let $\cat{Q}$ be a
small full sub-category of the category of comodules
$\lcomod{\coring{C}}$ whose underlying left $R$-modules are
finitely generated and projective. Denote by $\lambda_{\fk{q}}$
the coaction of $\fk{q}\in \cat{Q}$. Then one can directly apply
the above constructions, by putting $\chi(\fk{q}) = U(\fk{q})$,
where $U: \lcomod{\coring{C}} \rightarrow \lmod{R}$ is the left
forgetful functor. In this case, the left $\coring{C}$-coaction of
$\Sigma=\oplus_{\fk{q} \in \cat{Q}} \fk{q}$ is right $S$-linear,
while the right $\coring{C}$-coaction of ${}^{\vee}\Sigma$ is left
$S$-linear. Moreover, there is a canonical morphism of $R$-corings
defined by
\begin{equation}\label{can}
\xymatrix@R=0pt{\can{S}:
\Sigma \tensor{S}{}^{\vee}\Sigma \ar@{->}[r] & \coring{C} \\
u_{\fk{q}} \tensor{S} \varphi_{\fk{q}} \ar@{|->}[r] &
( \coring{C} \tensor{R} \varphi_{\fk{q}}) \,\circ \, \lambda_{\fk{q}}(u_{\fk{q}}). }
\end{equation}

Here $S$ is the induced ring from the category $\cat{Q}$, that is,
\begin{equation}\label{indring} S\,=\,
\oplus_{\fk{q},\,\fk{p} \in\cat{Q}\,}
\hom{\coring{C}}{\fk{q}}{\fk{p}}.
\end{equation}
However, the construction of the infinite comatrix coring, as well
as the canonical map $\can{}$, can be also performed for any
sub-ring of $S$ with the same set of orthogonal idempotents (i.e.
the $\fk{q}$'s identities).

Let us consider the $\Bbbk$-linear category $\Bbbk(\mathbb{N})$  whose objects are the natural numbers $\mathbb{N}$, and
homomorphisms sets are defined by
\begin{displaymath}
\hom{\Bbbk(\mathbb{N})}{n}{m} \,=\, \begin{cases}
         0,\text{ if }\, m \notin \{n,n+1\} \\ \Bbbk.1_{n}, \text{ if }\, n =m  \\  \Bbbk. \jmath_n^{n+1}, \text{ if }\, m=n+1,
        \end{cases}
\end{displaymath}
where the last two terms are free $\Bbbk$-modules of rank one. The
induced ring with enough orthogonal idempotents is the free
$\Bbbk$-module $B=\Bbbk^{(\mathbb{N})} \oplus
\Bbbk^{(\mathbb{N})}$ generated by the set
$\{\fk{h}_n,\fk{v}_n\}_{n \in \mathbb{N}}$, where $\fk{h}_n$ and
$\fk{v}_n$ corresponds to $1_n$ and $j_n^{n+1}$ respectively,
subject to the following relations:
\begin{eqnarray}
\fk{h}_n  \fk{h}_m &=& \delta_{n,\,m} \fk{h}_n,\,\, \forall m, n \in \mathbb{N} \qquad (\text{Kronecker delta}) \label{Eq:hv} \\ \nonumber
\fk{v}_n \fk{v}_m &=& \fk{v}_m \fk{v}_n \,\,=\,\,0,\,\, \forall m, n \in \mathbb{N}   \\ \nonumber
\fk{v}_n \fk{h}_{n+1} &=& \fk{v}_n \,\,=\,\, \fk{h}_n \fk{v}_n,\,\, \forall m, n \in \mathbb{N}.
\end{eqnarray}
In other words $B$ is the  ring of
$(\mathbb{N}\times \mathbb{N})$-matrices over $\Bbbk$ of the form
\begin{equation}\label{Ring-B}
\begin{pmatrix}  \Bbbk & \Bbbk & 0  & 0& & & \\ 0 & \Bbbk & \Bbbk
& 0 & &  & \\ 0& 0 & \Bbbk & \Bbbk & & &
\\  \vdots & & \ddots & \ddots & \ddots &  &
\\  & & &  & \Bbbk & \Bbbk &  \\  & & & & \ddots & \ddots & \ddots
 \end{pmatrix}
\end{equation}
consisting of matrices with only  possibly two non-zero entries in
each row: $(i,i)$ and $(i,i+1)$. It is clear that the category of
unital left $B$-modules is isomorphic to the category $\chain{\Bbbk}$ of
chain complexes of $\Bbbk$-modules. Precisely, this
isomorphism functor $\oO$ sends every chain complex $(V_{\bullet},
\partial^V)$ to its associated differential graded $\Bbbk$-module
$\oO(V_{\bullet}) \,=\,\oplus_{n \geq 0} V_n$ with the following
left $B$-action
$$ \fk{h}_{n} .\lr{ \sum_{n \geq 0}v_i} \,\,=\,\, v_n,  \quad \text{and }\, \fk{v}_n .\lr{\sum_{n \geq 0}v_i}\,\,=\,\, \partial^{V} (v_{n+1})$$
and acts in the obvious way on morphisms of chain complexes. The
inverse functor is clear.

By Proposition \ref{Q-comod},
we have a faithful functor $\chi: \Bbbk(\mathbb{N}) \to
\lcomod{\coring{\lL}}$ sending $n \to Q_n$, whose composition with
the left forgetful functor gives rise to a fiber functor $\chi:
\Bbbk(\mathbb{N}) \to add({}_RR)$.  Therefore, we can apply the above
process to construct an infinite comatrix $R$-coring
$Q\tensor{B}{}^{\vee}Q$ where $Q=\oplus_{n \in \mathbb{N}} Q_n$
and ${}^{\vee}Q=\oplus_{n \in \mathbb{N}} {}^*Q_n$ are given by
the cochain complex of Subsection \ref{cochain-Q}.

Since each of the $Q_n$'s has a structure of $R$-bimodule  for which the differential $\partial_{\bullet}$
is $R$-bilinear, we deduce that $Q\tensor{B}{}^{\vee}Q$ is an $\Re$-bimodule with actions
\begin{equation}\label{QQ}
 (r\tensor{}s^o) \, .\, (u_n\tensor{B}\varphi_n) \, . \,(p\tensor{}q^o) \,\,=\,\, (r u_np)\tensor{B} (s \varphi_n p),
\end{equation}
for every $p,q,r,s \in R$ and $u_n \in Q_n$ and $\varphi_n \in {}^*Q_n$. In view of this $\Re$-biaction,
the infinite comatrix $R$-coring has ${}_{\Re}(Q\tensor{B}{}^{\vee}Q)$
as its underlying $R$-bimodule.

Next we will construct an $\Re$-ring structure on the
$\Re$-bimodule $(Q\tensor{B}{}^{\vee}Q)$. Part of this construction needs the
following general Lemma which can be found, under a slightly
different form, in \cite{Deligne:1990}, \cite{Bruguieres:1994},
and \cite{P.H. Hai:2008}. We adopt the following general
notations: For any small $\Bbbk$-linear category $\cat{C}$, we
denote by ${\rm Funct}_{f}(\cat{C},add({}_RR))$ the category of
$\Bbbk$-linear faithful functors valued in $add({}_RR)$, i.e. of fiber functors on $\cat{C}$. For any object $\chi: \cat{C} \to
add({}_RR)$, we denote by $\cat{L}(\chi)$ the associated infinite
comatrix $R$-coring defined by the isomorphism of \eqref{R(S)}. Lastly, we
consider $\Sigma: {\rm Funct}_f(\cat{C},add({}_RR)) \to
\rmod{S(\cat{C})}$ the canonical functor to the category of right
unital $S(\cat{C})$-modules (recall that $S(\cat{C})$ is the
induced ring of $\cat{C}^{o}$).  That is,
\begin{equation}\label{Sigma}
\Sigma(\chi):=
\underset{\fk{c} \, \in \, \cat{C}}{\oplus} \fk{c}^{\chi}, \qquad \Sigma(\gamma):=  \underset{\fk{c} \, \in \, \cat{C}}{\oplus} \gamma_{\fk{c}}
\end{equation}
for every fiber functor $\chi$ and natural transformation $\gamma$
between fibred functors.

\begin{lemma}\label{omega-chi}
Let $\cat{A}$ be a small $\Bbbk$-linear category and let $\chi_1,
\chi_2: \cat{A} \to \bimod{R}$  be functors with images in
$add(_RR)$. Define $(\chi_1\tensor{R}\chi_2): \cat{A} \times
\cat{A} \to \bimod{R}$ by setting
$$(\chi_1\tensor{R}\chi_2)(\fk{p},\fk{q}) \,\,=\,\, \chi_1(\fk{p})
\tensor{R}\chi_2(\fk{q}), \text{ for } \fk{p}, \fk{q} \in \cat{A}.$$ Then
\begin{enumerate}[(i)]
 \item There is a left  $\Re$-linear isomorphism
$\cat{L}(\chi_1\tensor{R}\chi_2) \,\, \cong \,\, \cat{L}(\chi_1) \tensor{\Re} \cat{L}(\chi_2)$.
\item For every $R$-bimodule $M$, there is a natural isomorphism
\begin{footnotesize}
$$
\xymatrix@R=0pt{ {\rm Nat}\lr{(\chi_1\tensor{R}\chi_2),\, M\tensor{R}(\chi_1\tensor{R}\chi_2)} \ar@{->}[rr] & & \hom{R-R}{\cat{L}(\chi_1)
\tensor{\Re}\cat{L}(\chi_2)}{M} \\ \sigma \ar@{|->}[rr] & &
\Lr{(u\tensor{S}\varphi)\tensor{\Re} (v\tensor{S}\psi) \mapsto \sum_i m_i\varphi(p_i\psi(q_i)) }
}
$$
\end{footnotesize}
where $\sigma_{(\fk{p},\, \fk{q})} (u\tensor{R}v)\,=\, \sum_i m_i\tensor{R}p_i\tensor{R}q_i \in M\tensor{R}\fk{p}\tensor{R}\fk{q}$,
for every $u \in \fk{p}$, $\varphi \in {}^*\fk{p}$,  $v\in \fk{q}$, $\psi \in {}^*\fk{q}$ and
$(\fk{p},\fk{q}) \in \cat{A} \times \cat{A}$.
\end{enumerate}
\end{lemma}
\begin{proof}
Straightforward.
%
\end{proof}

Let us come back to our situation.  We are considering the functor
\begin{equation}\label{Eq:Chi}
\chi: \Bbbk(\mathbb{N}) \longrightarrow \lcomod{\coring{\lL}}, \text{ sending } \,n \longmapsto
Q_n.
\end{equation}
On the one hand, we already observed that the composition of
$\chi$ with the left forgetful functor gives rise to a fiber
functor $ \Bbbk(\mathbb{N}) \to add({}_RR)$. On the other hand, we
can consider also the fiber functor $\chi: \Bbbk(\mathbb{N}) \to
\bimod{R}$ obtained by composing the functor $\chi:
\Bbbk(\mathbb{N}) \to \lcomod{\coring{\lL}}$ with the functor
$\lcomod{\coring{\lL}}\rightarrow \bimod{R}$ stated in Lemma
\ref{Mono-comod-1}. Therefore, it is clear from
Lemma \ref{omega-chi}, that there is a bijective correspondence between  multiplications on $\cat{L}(\chi)\,=\, (Q\tensor{B}{}^{\vee}Q)$ and
natural transformations $(\chi\tensor{R}\chi) \to \cat{L}(\chi)
\tensor{R} (\chi\tensor{R}\chi)$. One of those natural transformations can be constructed
using the left $\cat{L}(\chi)$-coaction on the $Q_n$'s, as defined in
\eqref{td-lambda}. Thus we have the following statement

\begin{lemma}\label{natural-transf}
Let $Q_{\bullet}$ be the cochain complex of Subsection
\ref{cochain-Q}, and $(Q\tensor{B}{}^{\vee}Q)$ the
associated $R$-coring. Then there is a natural transformation
$(\chi\tensor{R}\chi) \to \cat{L}(\chi) \tensor{R}
(\chi\tensor{R}\chi)$ given by: $\td{\lambda}_{n,m}:
Q_n\tensor{R}Q_m \to (Q\tensor{B}{}^{\vee}Q) \tensor{R}
(Q_n\tensor{R} Q_m)$
\begin{footnotesize}
$$
u_n\tensor{R}u_m \longmapsto  \sum_{\alpha,\, \beta}
\Lr{(u_n\tensor{A}u_m)\tensor{B}({}^*\omega_{n,\,\alpha}\star \partial {}^*\omega_{m,\, \beta} )
+ (u_n\tensor{A}\partial u_m) \tensor{B} ({}^*\omega_{n,\,\alpha}\star {}^*\omega_{m,\, \beta} )}\,\tensor{R}\,
\lr{\omega_{n,\,\alpha}\tensor{R}\omega_{m\, \beta} }
$$
\end{footnotesize}
for every $n, m \geq 1$, and by $\td{\lambda}_{0,\, n}\,=\,
\td{\lambda}_{n,\, 0}: Q_n \to (Q\tensor{B}{}^{\vee}Q)
\tensor{R} Q_n$, sending $u_n \longmapsto \sum_{\alpha}
(u_n\tensor{B}{}^*\omega_{n,\, \alpha}) \tensor{R} \omega_{n,\,
\alpha}$, where  $\{(\omega_{n,\,\alpha},\, {}^*\omega_{n,\,
\alpha})\}$ is a dual basis for ${}_RQ_n$, $n \geq 1$.
\end{lemma}
\begin{proof}
This is a routine computation using definitions and dual bases notions.
\end{proof}
The following lemma will be  used in the sequel.

\begin{lemma}\label{suma-0}
Let $\{\omega_{n,\alpha},{}^*\omega_{n,\alpha})\}_{\alpha}$ be a dual basis for ${}_RQ_n$ with $n>0$. Then, for every element $u_n \in Q_n$,
$u_m \in Q_m$, and $\varphi_n \in {}^*Q_n$, $\varphi_m \in {}^*Q_m$, we have
\begin{footnotesize}
$$ \sum_{\alpha,\, \beta}
\Lr{(u_n\tensor{A}\partial u_m)\tensor{B}({}^*\omega_{n,\alpha}\star{}^*\omega_{m,\beta})} \times_R \Lr{(\omega_{n,\alpha}\tensor{A}
\omega_{m,\beta})\tensor{B}(\varphi_n\star\partial\varphi_m)} \,\,=\,\, 0
$$
\end{footnotesize}
and
\begin{footnotesize}
$$
\sum_{\alpha,\, \beta} \Lr{(u_n\tensor{A}u_m)\tensor{B}({}^*\omega_{n,\alpha}\star\partial{}^*\omega_{m,\beta})}
\times_R \Lr{(\omega_{n,\alpha}\tensor{A}\partial\omega_{m,\beta})\tensor{B}(\varphi_n\star\varphi_m)} \,\,=\,\,0
$$
\end{footnotesize}
as elements in the $\Re$-bimodule $\cross{(Q\tensor{B}{}^{\vee}Q)}{(Q\tensor{B}{}^{\vee}Q)}$.
\end{lemma}
\begin{proof}
Straightforward.
\end{proof}

We then arrive to the $\Re$-ring structure of $(Q\tensor{B}{}^{\vee}Q)$.

\begin{proposition}\label{D-algebra}
There is a structure of $\Re$-ring on $\dD:=(Q\tensor{B}{}^{\vee}Q)$ given by the extension of rings $\Re \to \dD$ sending
$r\tensor{}s^o \mapsto (r\tensor{B}s)$ (i.e. $\iota_0(r)\tensor{B}\iota_0(s)$, $\iota_0$ is the canonical injection),
where the multiplication of $\dD$ is defined by the following rules:
for every pair of generic elements $(u_n\tensor{B}\varphi_n)$ and $(u_m\tensor{B}\varphi_m)$ of $\dD$ with $n,m >0$, we set
$$
(u_n\tensor{B}\varphi_n) \,.\, (u_m\tensor{B}\varphi_m) \,=\,
\lr{(u_n\tensor{A}\partial u_m)\tensor{B}(\varphi_n\star\varphi_m) }
+ \lr{(u_n\tensor{A}u_m)\tensor{B}(\varphi_n\star\partial\varphi_m)}
$$
and
$$ (u_n\tensor{B}\varphi_n) \,.\, (r \tensor{B}s) \,\,=\,\, (u_n r\tensor{B} s \varphi_n), \quad
(r \tensor{B}s)\,.\,(u_n\tensor{B}\varphi_n)\,\,=\,\, (ru_n\tensor{B}\varphi_n s),\,\,\, \forall r, s \in R.$$
\end{proposition}
\begin{proof}
Using Lemmas \ref{dual-basis} and \ref{suma-0}, one can show that
each of the maps $\td{\lambda}_{n, m}$ given in Lemma
\ref{natural-transf} is coassociative with respect to the
comultiplication of $Q\tensor{B}{}^{\vee}Q$. Hence, its image by
the natural isomorphism  of Lemma \ref{omega-chi} leads to the
stated associative multiplication. The unitary property is clear.
\end{proof}

\begin{remark}
As we have seen, the construction of an $\Re$-ring structure on
$\dD$ is not an immediate task.  Part of this difficulty is clearly due to the fact that
the natural transformations which lead to the multiplications on $\dD$ are not easy to  construct. The other part is probably due to the fact
that, although the category $\Bbbk(\mathbb{N})$ is a monoidal
category, the fiber functor $\chi: \Bbbk(\mathbb{N}) \to
\bimod{R}$ given by the complex $Q_{\bullet}$ is not
strong monoidal since the local ''comultiplication'' maps
$Q_{n+m} \to Q_n\tensor{R}Q_m$, $m,n\geq 1$, see Remark
\ref{Q-normalization}, do not necessary form  a natural
isomorphisms. Thus $\chi$ does not satisfy the usual condition of a fibre functor, namely, being a strict monoidal functor.  Of course, this has prevented us from directly using
general results already existing in the literature, for example
\cite{P.H. Hai:2008}.
\end{remark}

\begin{proposition}\label{D-bialg}
Set $\dD:= {}_{\Re}(Q\tensor{B}{}^{\vee}Q){}_{\Re}$, where $Q_{\bullet}$ is the cochain complex defined in Subsection \ref{cochain-Q}.
Then $\dD$ has a structure of left $R$-bialgebroid.
\end{proposition}
\begin{proof}
Is a routine computation which uses Lemmas \ref{suma-0} and  \ref{dual-basis}, as well as Proposition \ref{D-algebra}.
\end{proof}

\subsection{The isomorphism between comatrices and coendomorphisms bialgebroids.}

Now, we come back to the canonical map. As  mentioned in the
preamble of the previous subsection, there is a canonical map
given explicitly by \eqref{can}. Thus, using the $\lL$-coactions
of Proposition \ref{Q-comod}, we have a morphism of $R$-corings
$\can{B}: \dD^l \longrightarrow \lL^l$ sending
\begin{footnotesize}
\begin{equation}\label{can-D}
(u_{n} \tensor{B} \varphi_{n}) \longmapsto \underset{i_0,\,i_1,\cdots,\,i_{n-1}}{\sum} \pi(a_0\tensor{} {}^*e_{i_0})\cdots
\pi(a_{n-1}\tensor{ } {}^*e_{i_{n-1}}) \varphi_n \lr{ e_{i_0}\partial e_{i_1} \tensor{A} \cdots \tensor{A}\partial e_{i_{n-1}} },
\end{equation}
\end{footnotesize}
where $u_n = a_0\partial a_1\tensor{A}\cdots\tensor{A}\partial a_{n-1} \in Q_n$, and $\can{B}(r\tensor{B}s)=\pi(r\tensor{}s^o)$, for
$r ,s \in R$.

Our next goal is to show that
$\can{B}$ is  an isomorphism of left $R$-bialgebroids.  To this end, we will need the following proposition.

\begin{proposition}\label{generadores}
For every $n \geq 1$, $u_n=a_0\partial a_1\tensor{A}\cdots\tensor{A}\partial a_{n-1} \in Q_n$ and  $\varphi_n \in {}^*Q_n$, we have
the following equality
\begin{footnotesize}
$$
(u_n\tensor{B}\varphi_n) \,=\,
\underset{i_0,i_1,\cdots,\,i_{n-1}}{\sum}
\Lr{(a_0\tensor{B}{}^*e_{i_0}).(a_1\tensor{B}{}^*e_{i_1}) \cdots
(a_{n-1}\tensor{B}{}^*e_{i_{n-1}}) }
\varphi_n(e_{i_0}\partial
e_{i_1}\tensor{A}\cdots\tensor{A}\partial e_{i_{n-1}})
$$
\end{footnotesize}
viewed as elements in the left $\Re$-module $\dD^l$. In particular, $\dD$ is generated, as an $\Re$-ring, by the image of $\Re$ and the set of elements
$\{(e_i\tensor{B}{}^*e_j)\}_{i,\,j}$ (recall that $\{(e_i,{}^*e_i)\}_i$ is a dual basis of ${}_RA$).
\end{proposition}
\begin{proof}
It follows by induction, using the dual basis of the $Q_n$'s given in Lemma \ref{dual-basis}.
\end{proof}

\begin{theorem}\label{can-bijective}
The canonical map $\can{B}: \dD \to \lL$ of \eqref{can-D} is an
isomorphism of left $R$-bialgebroids.
\end{theorem}
\begin{proof}
First we will show that $\can{B}$ is a multiplicative map. By
Proposition \ref{generadores} this is equivalent to show that
\begin{equation}\label{iguales}
\can{B}(a \tensor{B} \varphi)\,\can{B}(u_n\tensor{B}\varphi_n) \,=\,  \can{B}\lr{(a \tensor{B}\varphi)\,(u_n\tensor{B}\varphi_n)},
\end{equation}
for every $a \in A$, $\varphi \in {}^*A$, $u_n \in Q_n$,
$\varphi_n \in {}^*Q_n$ with $n \geq 1$.
Equality \eqref{iguales}, is  proved by direct computation. Since $\can{B}$ preserves the unit, we deduce
that $\can{B}$ is a morphism of $\Re$-rings.
The inverse of $\can{B}$ is constructed as follows. It is clear that
the map $\zeta:(A\tensor{}{}^*A)  \to \dD$ sending
$a\tensor{}\varphi \mapsto a\tensor{B}\varphi$ is an
$\Re$-bilinear map. Therefore, it is canonically extended to the
tensor algebra $\zeta: \tT_{\Re}((A\tensor{}{}^*A)) \to
\dD$, as $\dD$ is an $\Re$-ring. Now, for every $a, b \in A$ and
$\varphi \in {}^*A$, one shows that
$$
 \zeta\lr{\sum_i(a\tensor{}e_i\varphi)\tensor{\Re}(b\tensor{}{}^*e_i) } \,\,=\,\,
\zeta(ab\tensor{}\varphi),
$$
where $\{(e_i,{}^*e_i)\}_i$ is the dual basis of ${}_RA$.
This means that $\zeta$ factors throughout the canonical projection $\pi: \tT_{\Re}\lr{(A\tensor{}{}^*A)} \to \lL$,
and so we have an algebra map $\zeta: \lL \to \dD$. Given $a\in A$ and $\varphi \in {}^*A$,
we have
\begin{eqnarray*}
\can{B} \circ \zeta(\pi(a\tensor{}\varphi) )&=& \can{B}(a\tensor{B}\varphi) \,\,=\,\, \sum_i \pi(a\tensor{B}{}^*e_i)\varphi(e_i) \\
&=& \pi\lr{\sum_ia\tensor{}{}^*e_i\varphi(e_i)} \,\, =\,\,
\pi(a\tensor{}\varphi).
\end{eqnarray*}
This implies that $\can{B} \circ \zeta \,=\, id_{\lL}$. Now, take
$u_n \in Q_n$, $n \geq 1$, of the form $u_n=a_0\partial
a_1\tensor{A}\cdots\tensor{A}\partial a_{n-1}$ and $\varphi_n \in
{}^*Q_n$. Then,  by Proposition \ref{generadores}, we have
$$
\zeta \circ \can{B} (u_n\tensor{B}\varphi_n)\,\,=\,\, u_n\tensor{B}\varphi_n,
$$
and this shows that $\zeta \circ \can{B}\,=\,id_{\dD}$.
\end{proof}

\begin{corollary}\label{coendo}
Let $(\lL^l)^*$ be the right convolution ring of the $R$-coring $\lL^l$. Then there is an isomorphism of rings
$(\lL^l)^* \,\cong\, {\rm End}(Q_B)$.
\end{corollary}
\begin{proof}
We know that each $\fk{h}_n{}^{\vee}Q\,=\, {}^{*}Q_n$ is a finitely
generated and projective right $R$-module, where the $\fk{h}_i$'s are defined in \eqref{Eq:hv}.  The same property
holds true for each right $R$-module of the form
$\fk{e}_{i_1,i_n}{}^{\vee}Q$, where $\fk{e}_{i_1,i_n}=
\fk{h}_{i_1} + \cdots +  \fk{h}_{i_n}$. This means that the unital
bimodule ${}_B{}^{\vee}Q_R$ satisfies the second condition of
\cite[Proposition 5.1]{El Kaoutit:2006} for each idempotent which
belongs to the set of local units of $B$.  Therefore we have, as in
the proof of \cite[Proposition 5.1]{El Kaoutit:2006}, that the
functor $-\tensor{B}{}^{\vee}Q$ is left adjoint to $-\tensor{R}Q$.
Hence
$$\hom{-R}{\dD}{R}=\hom{-R}{Q\tensor{B}{}^{\vee}Q}{R}\cong \hom{-B}{Q}{Q}.$$ Now,
we conclude by Theorem \ref{can-bijective}.
\end{proof}

\section{Categories of comodules and chain complexes of modules.}\label{Main Results}
This section contains our main results, namely Theorems
\ref{equivalencia-B}, \ref{equivalencia-C} and \ref{Q_B-flat}. As
a consequence of these results, we obtain that the category of chain complexes of
left $R$-modules is always equivalent to the category of left
comodules over a quotient $R$-coring of the left $R$-bialgebroids
$\lL(A)$ constructed in Example \ref{main-example}. When $R$ is
commutative, this quotient inherits a left $R$-bialgebroid
structure from $\lL(A) $, and the stated equivalence is actually
a monoidal equivalence. This will clarify the equivalence of categories already constructed by Pareigis  and Tambara, \cite{Pareigis:81, Tambara}.

\smallskip

\subsection{Monoidal equivalence between chain complexes of $\Bbbk$-modules and left $\lL$-comodules.}\label{Result-1}
In this subsection we will use the isomorphism of bialgebroids stated in Theorem \ref{can-bijective}
to show that the following are equivalent: 1) $Q_B$ is faithfully flat, 2) the underlying module ${}_{R\tensor{}1^o}\lL$ of $\lL$
is flat and the functor $Q\tensor{B}- :\lmod{B} \to \lcomod{\lL}$ is a monoidal equivalence of categories. This is our first main result,
and stated below  as Theorem \ref{equivalencia-B}.

\begin{remark}\label{monoidal-B}
Let $B=\Bbbk^{(\mathbb{N})}\oplus\Bbbk^{(\mathbb{N})}$ be the ring
with enough orthogonal idempotents associated to the small
$\Bbbk$-linear category $\Bbbk(\mathbb{N})$ considered in
Subsection \ref{comatrix}, see \eqref{Ring-B}.  We have already observed in Subsection
\ref{comatrix} that the category of  unital left $B$-modules
$\lmod{B}$ is in a canonical way isomorphic to the category $\chain{\Bbbk}$ of
chain complexes of $\Bbbk$-modules. Therefore,
$\lmod{B}$ inherits a structure of monoidal category. Recall that
$B$ is generated as a free $\Bbbk$-module by the set of elements
$\{\fk{h}_n,\fk{v}_n\}_{n \in \mathbb{N}}$ with $\{\fk{h}_n\}_{n
\in\mathbb{N}}$ as a set of orthogonal idempotents given by \eqref{Eq:hv}. The
multiplication of two object $X, Y \in \lmod{B}$, is then given by
$$X \ominus Y\,\,=\,\, \underset{n \in\mathbb{N}}{\bigoplus} \lr{\underset{i+j=n}{\oplus} \fk{h}_iX\tensor{}\fk{h}_jX}, $$
That is, $\fk{h}_n (X\ominus Y)
\,\,=\,\,\oplus_{i+j=n}\fk{h}_iX\tensor{}\fk{h}_jY$, for every $n
\in \mathbb{N}$, and for every $k \geq 1$, $l \geq 1$ with
$k+l=m$, we have
$$ \fk{v}_{m-1} (\fk{h}_kx \tensor{}\fk{h}_ly)\,\,=\,\, \fk{v}_{ k-1}x \tensor{}\fk{h}_ly + (-1)^k \fk{h}_kx \tensor{}\fk{v}_{ l-1}y,$$
(i.e. the Leibniz rule), and $$\fk{v}_{ m-1} (\fk{h}_0 x\tensor{}\fk{h}_my)\,=\, \fk{h}_0 x\tensor{}\fk{v}_{ m-1}y,
\quad \fk{v}_{ n} (\fk{h}_{ n} x\tensor{}\fk{h}_0y)\,=\, \fk{v}_{ n-1} x\tensor{}\fk{h}_0y$$
for every $x \in X$, $y \in Y$,   and $m,n \geq 1$. The multiplication of $B$-linear maps is obvious.
The unit object is the left unital $B$-module $\Bbbk_{[0]}$ whose underlying $\Bbbk$-module is $\Bbbk$, and whose
$B$-action is given by
$$ \fk{h}_n \Bbbk_{[0]} \,\,=\,\, \begin{cases}
                                   0,\, \text{ if } n \neq 0 \\ \Bbbk,\,\text{ if } n=0.
                                  \end{cases}
$$
\end{remark}
We know that the cochain complex $Q_{\bullet}$ of Subsection \ref{cochain-Q} induces an $\lL$-comodule $Q\,=\,\oplus_{n \in \mathbb{N}} Q_n$
whose coaction is easily seen to be right $B$-linear. Thus, $Q\tensor{B}-:\lmod{B} \to \lcomod{\lL}$, acting in the obvious way,
is a well-defined functor. This functor is in fact  monoidal
\begin{lemma}\label{Q-monoidal}
Consider the monoidal categories $\lmod{B}$ and $\lcomod{\lL}$,
with structure, respectively, given in  Remark \ref{monoidal-B} and Lemma
\ref{Mono-comod-1}. Then $Q\tensor{B}-:\lmod{B} \to
\lcomod{\lL}$ is a monoidal functor, with structure
$$
\Gamma^2_{X,Y}: (Q\tensor{B}X) \tensor{R} (Q\tensor{B}Y)  \longrightarrow Q\tensor{B}(X\ominus Y),
$$
explicitly given by
\begin{multline*}
\Gamma^2_{X,Y}\lr{(u_n\tensor{B}\fk{h}_n x)\tensor{R}(u_m\tensor{B}\fk{h}_my)}  \\
\,=\,\begin{cases}
 (u_n\tensor{A}u_m) \tensor{B}(\fk{h}_n x \tensor{} \fk{v}_{m-1} y) \,+\, (u_n\tensor{A}\partial u_m)
\tensor{B}(\fk{h}_n x \tensor{} \fk{h}_{m} y),\,\,n,m \geq 1 \\ \,\, \\
u_nu_m \tensor{B}(\fk{h}_nx \tensor{}\fk{h}_m y), \quad n=0 \text{ or } m = 0,
\end{cases}
\end{multline*}
for every $u_n \in Q_n$, $u_m \in Q_m$, $x \in X$ and $y \in Y$,
and  $\Gamma^0: R \rightarrow Q\tensor{B}\Bbbk_{[0]}$ sending $r \mapsto r\tensor{B}\fk{h}_01$.
\end{lemma}
\begin{proof}
The fact that $\Gamma^2_{X,Y}$ is a well-defined map comes from
the observation that the right $R$-action of $Q\tensor{B}X$ as
left $\lL$-comodule is given by the right $R$-action of $Q$ viewed
as left $\lL$-comodule. That is, the one given by the rule
\eqref{right}. Now, it is easily seen that the right $R$-action of
$Q$ given by \eqref{right} is exactly the right $R$-action of $Q$
we started with (i.e. that which comes from  the inclusion
${}_RK_R \subset A\tensor{R}A$). A direct computation, using Lemma
\ref{Q-comod-1}, shows that $\Gamma^2_{X,Y}$ is left
$\lL$-colinear, for each $X,Y$. We leave to the reader the proof
of the associativity and unitary properties of the pair $(\Gamma^2_{-,-},
\Gamma^0)$.
\end{proof}

Our first main result is the following.

\begin{theorem}\label{equivalencia-B}
Let $R$ be an algebra over a commutative ground ring $\Bbbk$, and
$A$ an $R$-ring which is finitely generated and projective as left
$R$-module. Consider the associated left $R$-bialgebroid
given in Proposition \ref{pro: bialgebra} and let
$B=\Bbbk^{(\mathbb{N})} \oplus \Bbbk^{(\mathbb{N})}$ be the ring
with enough orthogonal idempotents of \eqref{Ring-B}. Consider the
cochain complex $Q_{\bullet}$ of Subsection \ref{cochain-Q} with
its canonical right unital $B$-action and left $\lL$-coaction.
Then the following statements are equivalent
\begin{enumerate}
 \item The right module $\lL^l_R$ is  flat and  the functor $Q\tensor{B}-: \lmod{B} \longrightarrow \lcomod{\lL}$
is an equivalence of monoidal categories;
\item $Q_B$ is a faithfully flat unital module.
\end{enumerate}
\end{theorem}
\begin{proof}
The monoidal condition is, by Lemma \ref{Q-monoidal}, always
satisfied, so it can be omitted in the proof.
Henceforth, we only need to show that $\lL^l_R$ is flat and
$Q\tensor{B}-$ is an equivalence, if and only if $Q_B$ is a
faithfully flat module. By the left version of the
generalized faithfully flat descent Theorem \cite[Theorem 5.9]{El
Kaoutit/Gomez:2004}, we know that $Q_B$ is faithfully flat if and
only if $\dD^l_{R}={}_{1\tensor{}R^o}(Q\tensor{B}{}^{\vee}Q)$ is
flat and $Q\tensor{B}-: \lmod{B} \to \lcomod{\dD^l}$ is an
equivalence of category. We then conclude by Theorem
\ref{can-bijective}.
\end{proof}

Notice that, when $Q_B$ is faithfully flat, the inverse of the
functor  of $Q\tensor{B}-:\lmod{B} \to \lcomod{\lL}$ is given by
the cotensor product ${}^{\vee}Q\square_{\lL}-: \lcomod{\lL} \to
\lmod{B}$. The structure of bicomodule on ${}^{\vee}Q$ is given as
follows.   Recall that $Q$ is in fact an $(\lL, B)$-bicomodule,
that is, the left $\lL$-coaction of $Q$ is right $B$-linear. So,
since each of the $Q_n$ , $n \geq 0$, is finitely generated and
projective left $R$-module, each of the left duals ${}^*Q_n$
admits a right $\lL$-coaction, for which
${}^{\vee}Q$  becomes a $(B,\lL)$-bicomodule.

The condition  $\lL^l_R$ is flat, stated in item $(1)$ of  Theorem
\ref{equivalencia-B}, seems to be redundant. But, although
we can deduce form the equivalence of categories that the category
of left $\lL$-comodule is abelian, we can not affirm that the
forgetful functor $\lcomod{\lL} \to \lmod{R}$ is left exact. Thus,
$\lL^l_R$ is not  necessarily a flat module, see
\cite[Proposition 2.1]{ElKaoutit/Gomez/Lobillo:2004c}.

Consider the category $\chain{\Bbbk}$ of  chain complexes of
$\Bbbk$-modules and  denote  by $\oO: \chain{\Bbbk} \to
\lmod{B}$ the canonical isomorphism of categories explicitly given in Subsection
\ref{comatrix}. In the case when $R=\Bbbk$ is a field, it is known
that $Q_B$ is always faithfully flat wherever $ {\rm
dim}_{\Bbbk}(A) < \infty$ (a complete proof for a non commutative
field, that is, a division ring is given in Theorem \ref{Q_B-flat}
below). We thus obtain the following corollary

\begin{corollary}{\cite[Theorem 4.4]{Tambara}}\label{Tambara-Thm}
Let $\Bbbk$ be a field and $A$ an $\Bbbk$-algebra such that $1<
{\rm dim}_{\Bbbk}(A) < \infty$. Consider the associated
coendomorphism $\Bbbk$-bialgebra $\lL$ given in Proposition
\ref{pro: bialgebra}. Then the category $\chain{\Bbbk}$ of chain
complexes of $\Bbbk$-modules  is monoidally
equivalent, via the functor $(Q\tensor{B}-) \circ \oO:
\chain{\Bbbk} \to \lcomod{\lL}$, to the category of left
$\lL$-comodules.
\end{corollary}
\begin{proof}
By the foregoing observations, this is a direct consequence of
Theorem \ref{equivalencia-B}.
\end{proof}

Explicitly, the composition of the functor given in Corollary
\ref{Tambara-Thm} with the forgetful functor $\lcomod{\lL} \to
\lmod{\Bbbk}$  is given as follows. For any chain complex $V_{\bullet}$ in
$\chain{\Bbbk}$, we have
$$Q\tensor{B}\oO(V_{\bullet}) \,\, =\,\,
\frac{\underset{n \geq 0}{\oplus} (Q_n \tensor{}V_n)}{\left\langle \partial u_n \tensor{} x_{n+1} -u_n\tensor{}\partial x_{n+1} \right\rangle}_{n \geq 0}.$$

\subsection{Equivalence between chain complexes of $R$-modules and $\bara{\lL}$-comodules.}\label{Result-2}
Our main aim here is to extend the result of Theorem
\ref{equivalencia-B} to the category $\chain{R}$ of chain
complexes over left $R$-modules. In other words, we are interested
in relating the category of chain complexes of left $R$-modules
and the category of left $\lL(A)$-comodules over the left
$R$-bialgebroid of Proposition \ref{pro: bialgebra}. Precisely, we
show an analogue of Theorem \ref{equivalencia-B} where $\lL$ is
replaced by its quotient $R$-coring $\bara{\lL}:
=\lL(A)/\left\langle 1_{\lL}(r\tensor{}1^o
-1\tensor{}r^o)\right\rangle_{r \,\in \, R} $ and the ring $B$ by
its extension  $C=R^{(\mathbb{N})}\oplus R^{(\mathbb{N})}$. This
is our second main result, see Theorem \ref{equivalencia-C}.  Of
course, in this case, the monoidal equivalence of categories is
reduced to an equivalence, unless
 the base ring $R$ is commutative and the
extension $A$ is an $R$-algebra.

Let $A$ be an $R$-ring and assume that ${}_RA$ is a finitely generated and projective module. Fix a dual basis $\{(e_i,\, {}^*e_i\}_i$ for ${}_RA$, and
consider $\lL:=\lL(A)$ the left $R$-bialgebroid of Proposition \ref{pro: bialgebra}
We denote by $\pi: \tT_{\Re}\lr{(A\tensor{}{}^*A)} \to \lL$
the canonical projection.

\begin{lemma}\label{J}
Let $\jJ$ be the left ideal of $\lL$ generated by the following set of elements
$$\left\{\underset{}{} \pi(a r\tensor{}\varphi) -\pi(a \tensor{} r \varphi)\right\}_{ a \in A,\, \varphi \in {}^*A,\,
r \in R.}$$
Then $\jJ$ is a coideal of the underlying $R$-coring $\lL^l$.
\end{lemma}
\begin{proof}
An easy computation shows that
$$\pi(a r\tensor{}\varphi)- \pi(a\tensor{} r \varphi) \,\,=\,\, \pi(a\tensor{}\varphi)\lr{r\tensor{}1^o \,-\, 1\tensor{}r^o},$$
for every elements $a \in A$, $\varphi \in {}^*A$ and $r \in R$.
Thus, $\jJ$ as left $\Re$-bimodule  is generated by the set
$\{\Sf{g}_r \,\,:=\,\, 1_{\lL}.(r\tensor{}1^o \,-\, 1\tensor{}r^o)
\}_{r \in R}$. For  arbitrary elements $x \in \lL$ and $r \in R$,
we get
$$
\varepsilon(x \Sf{g}_r)\,=\, \varepsilon \lr{ x. (1\tensor{}\varepsilon(\Sf{g}_r)^o)} \,=\, 0,
$$
as $\varepsilon(\Sf{g}_r)\,=\,0$. Hence, $\varepsilon(\jJ)\,=\,0$.
On the other hand, for every $r \in R$ , we have
$$\Delta(\Sf{g}_r)\,\,=\,\,  (1_{\lL}\tensor{R}1_{\lL}) (r\tensor{}1^o) - (1_{\lL}\tensor{R}1_{\lL})(1\tensor{}r^o).$$
Using these equalities we can show that, for every $x \in \lL$ and
$r \in R$, we have
\begin{eqnarray*}
\Delta(x \Sf{g}_r) &=& \sum_{(x)} x_{(1)}\tensor{R}x_{(2)}(r\tensor{}1^o) - \sum_{(x)}x_{(1)}\tensor{R}x_{(2)}(1\tensor{}r^o)  \\
&=& \sum_{(x)} x_{(1)}\tensor{R}x_{(2)}(r\tensor{}1^o \,-\, 1\tensor{}r^o) ,
\end{eqnarray*}
where $\Delta(x)\,=\, \sum_{(x)} x_{(1)}\tensor{R}x_{(2)}$.
Therefore, $(\bara{\pi}\tensor{R}\bara{\pi}) \circ \Delta (x
\Sf{g}_r)\,=\,0$, for every $x \in \lL$ and $r \in R$, where
$\bara{\pi}: \lL \to \lL /\jJ$ is the canonical projection. Thus
$\jJ$ is a  coideal of $\lL$.
\end{proof}

Denote by $\bara{\lL}:=\lL/\jJ$ the quotient $R$-coring and by
$\bara{\pi}: \lL \to \bara{\lL}$ the canonical projection. Notice
that $\bara{\pi}$ is also left $\lL$-colinear. Consider the
cochain complex $Q_{\bullet}$ of Subsection \ref{cochain-Q}. We
know, by Proposition \ref{Q-comod}, that each $Q_n$ is a left
$\lL$-comodule. Hence each of them is a left $\bara{\lL}$-comodule with coaction
$$\bara{\lambda}_n: Q_n \to \lL\tensor{R}Q_n \to \bara{\lL}\tensor{R}Q_n,\quad n\geq 0.$$

\begin{lemma}\label{Q-bicomod}
The $\lL$-coaction $\bara{\lambda}_n$ is
right $R$-linear. That is, $Q_n$ is an $(\bara{\lL},R)$-bicomodule
(here $R$ is considered as a the trivial $R$-coring).
\end{lemma}
\begin{proof}
For $n=0$ the statement is trivial since
$\bara{\lambda}_0(r)=(r\tensor{}1^o)\bara{\pi}(1_{\lL})=\bara{\pi}(1_{\lL}) (1 \tensor{}r^o)$, for every $r \in R$.
Take $n\geq 1$ and an element $u_n \in Q_n$ of the form $u_n=a_0\partial a_1 \tensor{A} \cdots \tensor{A}\partial a_{n-1} $.
Then, for every $r \in R$, we have

\begin{eqnarray*}
\bara{\lambda}_{n}( u_{n} r) &=& \sum_{\alpha,\, i_n} \bara{\pi}\lr{\pi(a_0\tensor{}{}^*e_{i_0})\cdots\pi(a_{n-1}r\tensor{}{}^*e_{i_{n-1}})}
\tensor{R}\omega_{n,\alpha},\\ & & \quad \text{ where } \alpha=(i_0,\cdots,i_{n-1}),\,\text{ and }\, \omega_{n,\alpha}=e_{i_0}\partial
e_{i_1} \tensor{A} \cdots \tensor{A}\partial e_{i_{n-1}} \\
&=& \sum_{\alpha} \bara{\pi}\lr{\pi(a_0\tensor{}{}^*e_{i_0})\cdots\pi(a_{n-1}\tensor{}r{}^*e_{i_{n-1}})}
\tensor{R}\omega_{n,\alpha} \\
&=& \sum_{\alpha}
(\bara{\pi}\tensor{R}Q_n)\Lr{\pi(a_0\tensor{}{}^*e_{i_0})\cdots\pi(a_{n-1}\tensor{}{}^*e_{i_{n-1}})\,
(1\tensor{}r^o) \,
\tensor{R}\,\omega_{n,\alpha} } \\
&=& \sum_{\alpha}
(\bara{\pi}\tensor{R}Q_n)\Lr{\pi(a_0\tensor{}{}^*e_{i_0})\cdots\pi(a_{n-1}\tensor{}{}^*e_{i_{n-1}})
\,
\tensor{R}\,\omega_{n,\alpha}r } \\
&=& \sum_{\alpha} \bara{\pi}\lr{\pi(a_0\tensor{}{}^*e_{i_0})\cdots\pi(a_{n-1}\tensor{}{}^*e_{i_{n-1}})}
\tensor{R}\omega_{n,\alpha}r \\
&=& \bara{\lambda}_n(u_n) r,
\end{eqnarray*}
where in the fourth equality we have used that each $Q_n$ is in fact a left $\times_R$-$\lL$-comodule, that is, $\lambda_n(Q_n) \subseteq \lL \times_R Q_n$. We then conclude by linearity.
\end{proof}

\begin{remark}\label{R-commutative}
The quotient $R$-coring $\bara{\lL}$ does not admit, in general, a
structure of left $R$-bialgebroid. However, if we assume that $R$
is commutative (i.e. a commutative $\Bbbk$-algebra) and that $A$
is an $R$-algebra, then the left ideal $\jJ$ is in fact a
two-sided ideal, since in this case we have the following
equalities
$$ \Sf{g}_r \pi(a\tensor{}\varphi) \,=\, \pi(a\tensor{}\varphi) \Sf{g}_r, \text{ for every }\, r \in  R,\, a \in A, \text{ and } \varphi \in A^*.$$
In view of this, a direct verification shows that $\bara{\lL}$ is an $R$-bialgebroid such that the canonical surjection $\bara{\pi}: \lL \to \bara{\lL}$
is a morphism of $R$-bialgebroids. Notice, that here the prefix ''left'' was removed before bialgebroid. This is due to the fact that $\lL$ is actually
an $(R\tensor{}R)$-algebra, that is, there is only one structure of $(R\tensor{}R)$-module.
\end{remark}

Let us consider the $\Bbbk$-linear category $R(\mathbb{N})$  whose
objects are the natural numbers $\mathbb{N}$ and homomorphisms
sets are defined by
\begin{equation}\label{RN}
\hom{R(\mathbb{N})}{n}{m} \,=\, \begin{cases}
         0,\text{ if }\, m \notin \{n,n+1\} \\ R.1_{n}\,=\, 1_{n}.R, \text{ if }\, n =m  \\  R. \jmath_n^{n+1}\,=\, \jmath_{n}^{n+1}. R ,
\text{ if }\, m=n+1.
        \end{cases}
\end{equation}
The last two terms are copies of ${}_RR_R$ viewed as an
$R$-bimodule which is free as left and right $R$-module of rank
one, generated by an invariant element. The composition is defined
using the regular $R$-biactions of ${}_RR_R$. The induced ring
with enough orthogonal idempotents is the free left $R$-module
$C=R^{(\mathbb{N})} \oplus R^{(\mathbb{N})}$ generated by elements
$\{\fk{h}_n,\fk{u}_n\}_{n \in \mathbb{N}}$ subject to the
following relations:
\begin{eqnarray}
\label{Eq:hu} \fk{h}_n  \fk{h}_m &=& \delta_{n,\,m} \fk{h}_n, \,\, \forall n,\,m \in \mathbb{N} \qquad (\text{Kronecker delta}) \\ \nonumber
\fk{u}_n \fk{u}_m &=& \fk{u}_m \fk{u}_n \,\,=\,\,0, \,\, \forall n,\,m \in \mathbb{N} \\
\nonumber \fk{u}_n \fk{h}_{n+1} &=& \fk{u}_n \,\,=\,\, \fk{h}_n \fk{u}_n, \,\, \forall n,\,m \in \mathbb{N}.
\end{eqnarray}
In other words $C$ is the ring of $(\mathbb{N}\times \mathbb{N})$-matrices over $R$ of the form
\begin{equation}\label{Ring-C}
C\,\,=\,\,
\begin{pmatrix}  R & R & 0  & 0& & & \\ 0 & R & R
& 0 & &  & \\ 0& 0 & R & R & & &
\\  \vdots & & \ddots & \ddots & \ddots &  &
\\  & & & 0 & R & R &  \\  & & & & \ddots & \ddots & \ddots
 \end{pmatrix}
\end{equation}
i.e. with possibly non-zero entries in each row: $(i,i)$ and
$(i,i+1)$. $C$ is also free as right $R$-module, since the
generators are invariant. One can easily check that the category
of chain complexes of left $R$-modules $\chain{R}$ is equivalent
to the category of unital left $C$-modules. Let $B$ be the ring
with enough orthogonal idempotents of \eqref{Ring-B}. There is a
morphism of rings $B  \to C$ with the same set of orthogonal
idempotents. In this way, we have by \cite[page 733]{El
Kaoutit:2006} the usual adjunction between left unital $B$-modules
and $C$-modules using restriction of scalars and the tensor
product functor $C\tensor{B}-$.

By Lemma \ref{Q-bicomod}, we have a morphism of rings $R \to {\rm
End}_{\bara{\lL}}(Q_n)$, for every $n\geq 0$.  This leads to a
faithful functor from the category $R(\mathbb{N})$ to the category
of $(\bara{\lL},R)$-bicomodules (here $R$ is considered as a
trivial $R$-coring) $\chi': R(\mathbb{N}) \to
{}_{\bara{\lL}}\Sf{Comod}_{R}$. The composition of $\chi'$ with
the forgetful functor gives rise then to a fiber functor $\omega:
R(\mathbb{N}) \to \bimod{R}$ whose image is in $add({}_RR)$.
Therefore, we can apply the constructions performed in Subsection
\ref{comatrix}. Thus, we have an infinite comatrix $R$-coring
$Q\tensor{C}{}^{\vee}Q$ together with a canonical map $\can{C}:
Q\tensor{C}{}^{\vee}Q \longrightarrow \bara{\lL}$ sending
\begin{footnotesize}
\begin{equation}\label{can-C}
\xymatrix@C=40pt {u_n\tensor{C}\varphi_n \ar@{->}^-{\can{C}}[r] & \sum_{i_0,\,\cdots,\, i_{n-1}}
\bara{\pi}\lr{\pi(a_0\tensor{}{}^*e_{i_0})\cdots\pi(a_{n-1}\tensor{}{}^*e_{i_{n-1}})} \varphi\lr{e_{i_0}\partial
e_{i_1} \tensor{A} \cdots \tensor{A}\partial e_{i_{n-1}}}. }
\end{equation}
\end{footnotesize}
Clearly we have a surjective map $\vartheta: Q\tensor{B}{}^{\vee}Q \to Q\tensor{C}{}^{\vee}Q$. Moreover,
we have a commutative diagram with exact rows relating the two $R$-corings morphisms $\can{B}$ and $\can{C}$
(see equations \eqref{can-D} and \eqref{can-C})
\begin{footnotesize}
\begin{equation}\label{can-diagram}
\xymatrix@C=60pt{ 0 \ar@{->}[r] & {\rm Ker}(\vartheta) \ar@{->}[r]
\ar@{-->}[dd] & Q\tensor{B}{}^{\vee}Q \ar@{->}^-{\vartheta}[r]
\ar@{->}^-{\can{B}}[dd] & Q\tensor{C}{}^{\vee}Q
\ar@{->}^-{\can{C}}[dd] \ar@{->}[r] & 0 \\ & & & &
\\ 0 \ar@{->}[r] & \jJ \ar@{->}[r] & \lL \ar@{->}^-{\bara{\pi}}[r]  & \bara{\lL} \ar@{->}[r] & 0 }
\end{equation}
\end{footnotesize}

\begin{proposition}\label{can_C}
In diagram \eqref{can-diagram}, we have the following equality  $\can{B}({\rm Ker}(\vartheta)) = \jJ$. In particular, the map
$\can{C}$ of equation \eqref{can-C} is an isomorphism of $R$-corings.
\end{proposition}
\begin{proof}
The inclusion $\can{B}({\rm Ker}(\vartheta)) \subseteq \jJ$ is
clear from the commutative diagram \eqref{can-diagram}.
Conversely, consider arbitrary elements $y \in \lL$ and $r \in R$.
We need to show that $y \Sf{g}_r \in \can{B}({\rm
Ker}(\vartheta))$, where $\Sf{g}_r$ are as in the proof of Lemma \ref{J}. There is no loss of generality if we assume
that $y=x \pi(a\tensor{}\varphi)$, for some $x \in \lL$ and $a \in
A$, $\varphi \in {}^*A$. Since $\can{B}$ is, by Theorem
\ref{can-bijective}, bijective, there exists $u \in
Q\tensor{B}{}^{\vee}Q$ such that $x=\can{B}(u)$.  In view of this,
$y \Sf{g}_r \,=\, \can{B}(u(ar\tensor{B}{
\varphi}-a\tensor{B}r\varphi))$, as $\can{B}$ is multiplicative.
We need to check that $\vartheta(u\,(ar\tensor{B}{
\varphi}-a\tensor{B}r\varphi))\,=\,0$. However,  this is directly obtained
from the following equality
$$\vartheta\lr{(u_n\tensor{B}\varphi_n)\,(ar\tensor{B}-a\tensor{B}r\varphi)}\,=\,0,\, \text{ for every  } u_n \in Q_n,\, \varphi_n \in {}^*Q_n,$$
whose proof follows by induction on $n$.
%
The last statement to prove is a consequence
of the first one, since the diagram \eqref{can-diagram} has exact
rows.
\end{proof}

Our second main result is the following

\begin{theorem}\label{equivalencia-C}
Let $R$ be an algebra over a commutative ground ring $\Bbbk$, and $A$ an $R$-ring which is finitely generated and projective as left $R$-module.
Consider the associated left $R$-bialgebroid $\lL$ stated in
Proposition \ref{pro: bialgebra} and $\jJ$ the coideal of $\lL$ generated by the set of elements $\{1_{\lL}(r\tensor{}1^o - 1\tensor{}r^o)\}_{r \in R}$.
Denote by $\bara{\lL}=\lL/\jJ$ the corresponding quotient $R$-coring.
Let $C=R^{(\mathbb{N})} \oplus R^{(\mathbb{N})}$ be the ring with enough orthogonal idempotents induced from the small $\Bbbk$-linear category
$R(\mathbb{N})$ defined by relations \eqref{RN}. Consider the cochain complex $Q_{\bullet}$ given in Subsection \ref{cochain-Q}
with its canonical right unital $C$-action and left $\bara{\lL}$-coaction as in Lemma \ref{Q-bicomod}.
Then the following statements are equivalent
\begin{enumerate}
 \item The right module $\bara{\lL}^l_R$ is  flat and  the functor $Q\tensor{C}-: \lmod{C} \longrightarrow \lcomod{\bara{\lL}}$
is an equivalence of categories;
\item $Q_C$ is a faithfully flat unital module.
\end{enumerate}
\end{theorem}
\begin{proof}
By the left version of the generalized faithfully flat descent
Theorem \cite[Theorem 5.9]{El Kaoutit/Gomez:2004}, we know that
$(Q\tensor{C}{}^{\vee}Q)_R$ is flat and $Q\tensor{C}-:\lmod{C} \to
\lcomod{Q\tensor{C}{}^{\vee}Q}$ is an equivalence of categories,
if and only if $Q_C$ is faithfully flat. We then deduced the
stated equivalence by using the isomorphism of $R$-corings
$\can{C}: Q\tensor{C}{}^{\vee}Q \cong \bara{\lL}$ established in
Proposition \ref{can_C}.
\end{proof}

Notice that, if $Q_C$ is faithfully flat, then the inverse functor of $Q\tensor{C}-: \lmod{C} \to \lcomod{\bara{\lL}}$
is given by the cotensor product ${}^{\vee}Q\square_{\bara{\lL}}-: \lcomod{\bara{\lL}} \to \lmod{C}$. Here the structure of
$(C,\bara{\lL})$-bicomodule of ${}^{\vee}Q$ is deduced, as was observed in Subsection \ref{Result-1}, from that of $Q$ using the fact that each of the $Q_n$'s is finitely generated
and projective left $R$-module.

\subsection{Conditions under which  $Q_{C}$ is faithfully flat.}\label{Result-3}
As was seen in Theorems \ref{equivalencia-B} and
\ref{equivalencia-C}, a sufficient and necessary condition for establishing an
equivalence of categories of left comodules and chain
complexes, is the faithfully flatness of the unital right module
$Q$. The proof of this fact is actually the most difficult task in
this theory. In this subsection we will analyze assumptions under
which $Q_C$ is faithfully flat.

The following is our third main result.
\begin{theorem}\label{Q_B-flat}
The notations and assumptions are that of Theorem
\ref{equivalencia-C}. Assume further that $A_R$ is finitely
generated and projective, and  the cochain complex $Q_{\bullet}$
is exact and splits, in the sense that, for every $m \geq 1$,
$Q_m\,=\, \partial Q_{m-1} \oplus \bara{Q}_{m}\,=\,{\rm
Ker}(\partial) \oplus \bara{Q}_m$ as right $R$-modules, for some
right $R$-module $\bara{Q}_m$.  Then $Q_C$ is a flat module.
Furthermore, if $\Bbbk$ is a field and $R$ is a division
$\Bbbk$-algebra, then $Q_C$ is faithfully flat.
\end{theorem}
\begin{proof}
We first consider the following family of right $R$-modules
$$
Q^{(m)} \,\,=\, \, \begin{cases} \partial Q_m \,\oplus\, \bara{Q}_m, \quad \text{ for } m \geq 1
\\ \, \\ \partial Q_0 \,\oplus \, Q_0, \quad \text{ for } m = 0
                                           \end{cases}
$$
which we claim to be a family of right unital flat $C$-modules.
Using this claim we can easily deduce that $Q_C$ is a flat module
since  we know that $Q_C \,=\, \oplus_{m \geq 0}\, Q^{(m)}_C$. The
structure of unital right $C$-module of each $Q^{(m)}$ is given as
follows: Denote by $\Sf{i}_m: \partial Q_m \to Q^{(m)}$,
$\bara{\Sf{i}}_m: \bara{Q}_m \to Q^{(m)}$ the canonical injections
and by $\Sf{j}_m$, $\bara{\Sf{j}}_m$ their canonical projections.
For every element $u^{(m)} \in Q^{(m)}$, we set
$$
u^{(m)}\, \fk{h}_n \,=\, \begin{cases} 0,\text{ if } n \notin \{m,m+1\} \\ \bara{\Sf{i}}_m\bara{\Sf{j}}_m(u^{(m)}), \text{ if }
 n=m \\ \Sf{i}_m\Sf{j}_m(u^{(m)}), \text{ if } n=m+1  \end{cases}\qquad
u^{(m)}\, \fk{u}_n \,=\, \begin{cases} 0,\text{ if } n \neq m \\ \Sf{i}_m\lr{\gamma_m\bara{\Sf{j}}_m(u^{(m)})}, \text{ if }
 n=m   \end{cases}
$$
where $\gamma_m: \bara{Q}_m \to Q_m \to \partial Q_m$. That is,  the obtained cochain complexes have the following form
\begin{equation}\label{Eq:Qm}
\xymatrix{Q^{(m)}_{\bullet}: 0 \ar@{->}[r] & 0\cdots \cdots 0\ar@{->}[r] & \bara{Q}_m \ar@{->}[dr] \ar@{-->}^-{\gamma_m}[rr] & & \partial Q_m \ar@{->}[r]
& 0 \ar@{->}[r] & 0\cdots \cdots
\\  &  &  & Q_m \ar@{->}[ur] &   &  & }
\end{equation}
Put $\fk{e}_{n,\,n+1}=\fk{h}_n + \fk{h}_{n+1}$, for every $n \geq 0$. These are idempotents elements in $C$, and the induced rings, i.e.
$\fk{e}_{n,\, n+1} C \fk{e}_{n,\, n+1}$ are all isomorphic to the upper-triangular matrices over $R$. That is, of the form

$$C_{n,\, n+1}:\,= \,\fk{e}_{n,\, n+1} C \fk{e}_{n,\, n+1} \,\,=\,\,
\begin{pmatrix}
R  & R \\
0  & R
                                                  \end{pmatrix}, \quad \text{ for every }\, n \in \mathbb{N}.
$$
It is clear that, for every $m \geq 0$, we have  $Q^{(m)}\fk{e}_{m,\, m+1} \,=\, Q^{(m)}$.
Therefore, there is an  isomorphism of right unital $C$-modules
\begin{equation}\label{isos-m}
 Q^{(m)}\fk{e}_{m,\, m+1} \underset{C_{m,\, m+1}}{\otimes} \fk{e}_{m,\, m+1}C \,\, \cong \,\, Q^{(m)}.
\end{equation}

Next we will show that each of the right $C_{m,\, m+1}$-modules
$Q^{(m)}\fk{e}_{m,\, m+1}=Q^{(m)} $ is finitely generated and
projective. This fact, combined with the isomorphisms
\eqref{isos-m}, establish the above claim.

For $m=0$, it is clear that the right $C_{0,\, 1}$-module
$$Q^{(0)} = R\oplus R = {\footnotesize\begin{pmatrix}  R & R \\ 0
& 0
\end{pmatrix} }={\footnotesize\begin{pmatrix}  1_R & 0 \\ 0 & 0
\end{pmatrix} }C_{0,\,1}$$ is finitely generated and projective. Now take $m \geq
1$, under the hypothesis $A_R$ is finitely generated and
projective, we can show,  as in Lemma \ref{dual-basis}, that each
right $R$-module $Q_m$ is also finitely generated and projective.
Thus, we can consider a dual basis $\{(\bara{q}_{m,k},
\bara{q}_{m,k}^*)\}_{k}$ for each right $R$-module $\bara{Q}_m$.
In this way, we have a right $C_{m,\, m+1}$-linear map
$$
\theta^*_{m,\,k}: Q^{(m)} \longrightarrow  C_{m,\, m+1} , \qquad \Lr{
 u^{(m)} \longmapsto  \begin{pmatrix} \bara{q}^*_{m,k}(\bara{\Sf{j}}_m(u^{(m)}) & \bara{q}^*_{m,k}(\bara{x}_m) \\ & \\ 0 & 0 \end{pmatrix} },
$$
where $\bara{x}_m \in \bara{Q}_m$ is the projection of $x_m \in
Q_m=\partial Q_{m-1} \oplus \bara{Q}_m$, defined by
$\Sf{j}_m(u^{(m)}) \,=\, \partial x_m \in \partial Q_m$. We should
mention that, under our assumptions,  the maps $\theta^*_{m,k}$
are well defined. Effectively, if there is some other element $y_m
\in Q_m$ such that $\Sf{j}_m(u^{(m)}) \,=\, \partial x_m \,=\,
\partial y_m$, then $x_m - y_m \in {\rm Ker}(\partial_{m}) \,=\,
\partial Q_{m-1}$ which means that they have equal image
$\bara{x}_m \,=\, \bara{y}_m$ in $\bara{Q}_m \cong Q_m/\partial
Q_{m-1}$. It is convenient to check that $\theta_{m,\, k}^*$ are
right $C_{m,\, m+1}$-linear. But first we will identify the right
module $\bara{Q}_m$ with the quotient of $Q_m$, $\bara{Q_m}
=Q_m/\partial Q_{m-1}$. The right $C_{m,\,m+1}$-action of
$Q^{(m)}$ is given as follows: Take an element $u^{(m)} \in
Q^{(m)}$ and write it in the form $u^{(m)}\,=\, (\bara{q_m},
\partial p_m)$ for some elements $q_m ,\, p_m \in Q_m$. Here
$\Sf{j}_m(u^{(m)})\,=\, \partial p_m$ and
$\bara{\Sf{j}}_m(u^{(m)})\,=\, \bara{q_m}$. So
$$
(\bara{q_m},\, \partial p_m)\, \begin{pmatrix} r_{11} & r_{12} \\ 0 & r_{22}\end{pmatrix} \,\,=\,\, \lr{\bara{q_m}\, r_{11},\, \partial q_m \, r_{12} +
\partial p_m r_{22} },
$$
for every element $\begin{pmatrix} r_{11} & r_{12} \\ 0 & r_{22}\end{pmatrix}$ in $C_{m,\, m+1}$. Therefore,
\begin{eqnarray*}
\theta_{m,\, k}^*\lr{(\bara{q_m},\, \partial p_m)\, \begin{pmatrix} r_{11} & r_{12}\\ &  \\ 0 & r_{22}\end{pmatrix}} &=&
\theta_{m,\, k}^* \lr{\bara{q_m}\, r_{11},\, \partial q_m \, r_{12} +
\partial p_m r_{22} } \\
&=& \begin{pmatrix} \bara{q}_{m,\,k}^*(\bara{q_m} \, r_{11}) &
\bara{q}_{m,k}^*\lr{\bara{q_m} \, r_{12} + \bara{p}_m \, r_{22}}\\
& \\ 0 & 0
\end{pmatrix} \\
&=& \begin{pmatrix} \bara{q}_{m,\,k}^*(\bara{q_m} \, r_{11}) &
\bara{q}_{m,k}^*(\bara{q_m} \, r_{12}) +
\bara{q}_{m,k}^*(\bara{p}_m \, r_{22}) \\ & \\ 0 & 0
\end{pmatrix} \\
&=& \begin{pmatrix} \bara{q}_{m,\,k}^*(\bara{q_m}) \,r_{11}  &
\bara{q}_{m,k}^*(\bara{q_m}) \, r_{12} +
\bara{q}_{m,k}^*(\bara{p}_m) \, r_{22} \\ & \\ 0 & 0
\end{pmatrix}  \\
&=& \begin{pmatrix} \bara{q}_{m,k}^*(\bara{q_m}) &
\bara{q}_{m,k}^*(\bara{p}_m ) \\ & \\ 0 & 0
\end{pmatrix} \, \begin{pmatrix} r_{11} & r_{12} \\ & \\ 0 & r_{22}\end{pmatrix} \\
&=& \theta_{m,\, k}^*(\bara{q_m},\, \partial p_m)\, \begin{pmatrix} r_{11} & r_{12}\\ &  \\ 0 & r_{22}\end{pmatrix}.
\end{eqnarray*}
Take an arbitrary element $(\bara{q_m}, \partial p_m) \in Q^{(m)}$, we have
\begin{eqnarray*}
(\bara{q_m},\, \partial p_m) &=& (\bara{q_m},\, 0) \,+ \, (0,\,\partial p_m) \\
&=& (\bara{q_m},\, 0) \,+ \, (\bara{ p_m},\,0)\, \begin{pmatrix} 0 & 1 \\ 0 & 0\end{pmatrix} \\
&=& \sum_k (\bara{q}_{m,\, k} \bara{q}_{m,\, k}^*(\bara{q_m}),\,0) \,+\,
\sum_k (\bara{q}_{m,\, k} \bara{q}_{m,\, k}^*(\bara{p_m}),\,0)\, \begin{pmatrix} 0 & 1 \\ 0 & 0\end{pmatrix} \\
&=& \sum_k (\bara{q}_{m,\, k},\,0)  \begin{pmatrix} \bara{q}_{m,\, k}^*(\bara{q_m})  & 0 \\ 0 & 0\end{pmatrix}    \,+\,
\sum_k (\bara{q}_{m,\, k} ,\,0) \, \begin{pmatrix} \bara{q}_{m,\, k}^*(\bara{p_m})  & 0 \\ 0 & 0\end{pmatrix}
\, \begin{pmatrix} 0 & 1 \\ 0 & 0\end{pmatrix} \\
&=& \sum_k (\bara{q}_{m,\, k},\,0)  \begin{pmatrix} \bara{q}_{m,\, k}^*(\bara{q_m})  & 0 \\ 0 & 0\end{pmatrix}    \,+\,
\sum_k (\bara{q}_{m,\, k} ,\,0) \, \begin{pmatrix} 0 & \bara{q}_{m,\, k}^*(\bara{p_m})   \\ 0 & 0\end{pmatrix} \\
&=& \sum_k (\bara{q}_{m,\, k},\,0)  \begin{pmatrix} \bara{q}_{m,\, k}^*(\bara{q_m})  & \bara{q}_{m,\, k}^*(\bara{p_m}) \\ 0 & 0\end{pmatrix} \\
&=& \sum_k (\bara{q}_{m,\, k},\,0) \, \theta_{m,\, k}^*\lr{\bara{q_m},\, \partial p_m},
\end{eqnarray*}
which shows that $\left\{\lr{(\bara{q}_{m,\, k},\, 0),\,\theta_{m,\, k}^*} \right\}_k$ is a dual basis for the right $C_{m,\, m+1}$-module $Q^{(m)}$, and
this finishes the proof of the main statement.

If we assume now that $\Bbbk$ is a field and $R$ is a division
$\Bbbk$-algebra, then one can show as follows that each
$Q^{(m)}\fk{e}_{m,\, m+1}$ is a progenerator in the category of
right $C_{m,\, m+1}$-modules. This will imply that
$Q^{(m)}\fk{e}_{m,\, m+1}\tensor{C_{m,\, m+1}}-: \lmod{C_{m,\, m+1}} \to \lmod{R}$ is a faithful
functor. By identifying each ring $C_{m,\, m+1}$ with the
upper triangular matrix ring $T:=\begin{pmatrix} R & R\\ 0 &
R\end{pmatrix}$, we know that $T\,=\,eT \oplus (1-e)T$, where $e$
is the obvious idempotent element. The structure of right
$T$-module of $Q^{(m)}$ is given by the decomposition $Q^{(m)}_T
=\partial Q_{m} \oplus \bara{Q}_m$ with a surjective canonical map
$\gamma_m: \bara{Q}_m \to \partial Q_m$ of \eqref{Eq:Qm}.  Since $R$ is a division
ring and each component of $Q^{(m)}$ is by assumption finite
dimensional with $d={\rm dim}_R(\bara{Q}_m) \leq {\rm
dim}_R(\partial Q_m)=d'$, we can split $Q^{(m)}$ as
$$Q^{(m)} \,\cong\, (eT)^{d} \oplus \lr{(1-e)T}^{d'-d},$$ and this shows that $Q^{(m)}_T$ is a progenerator.
Let $f: X \to Y$ be a morphism of right unital $C$-modules such
that $Q\tensor{C}f\,=\, 0$. Hence $Q^{(m)} \tensor{C}f\,=\, 0$,
for every $m\geq 0$, as $Q_C\,=\, \oplus_{m \geq 0} Q^{(m)}$.
Therefore, we have
$$0=Q^{(m)} \tensor{C} f \, \cong \,Q^{(m)}\fk{e}_{m,\, m+1} \underset{C_{m,\,m+1}}{\otimes} \fk{e}_{m,\,m+1}C \tensor{C} f, \, \,\, \forall m \geq 0
\,\,\Longrightarrow \,\,
\fk{e}_{m,\,m+1}C \tensor{C} f \, =\, 0, \,\, \forall m \geq 0. $$
This means that $\fk{h}_m C\tensor{C} f \,=\, 0$, for every $m \geq 0$, and so $f\,=\,0$. This shows that $Q\tensor{C}-$ is a faithful functor, which
completes the proof.
\end{proof}

\begin{remark}
As one can see, the hypothesis on the complex $Q_{\bullet}$ in
Theorem \ref{Q_B-flat}, is not easy to check. However, under
further conditions on the ring extension $R \to A$, this
hypothesis is satisfied. For instance, it is clear from Lemma
\ref{mu-splits} and Remark \ref{Remark-mu-splits} that it is
satisfied by assuming that the ring extension $R \to A$  splits
either in the category of right or left $R$-modules. Obviously
this includes the case when $A$ is free as right (or left)
$R$-module with $1_A$ as an element of the canonical basis. In
particular, this is the case when $R$ is a division ring.
\end{remark}

\begin{corollary}
Let $D$ be a division $\Bbbk$-algebra over a field $\Bbbk$, and
$A$ a $D$-ring which is finite dimensional as left and right
$D$-vector space with dimension $ \geq 2$. Consider the associated
left $D$-bialgebroid $\lL$  given by Proposition \ref{pro:
bialgebra} and its coideal $\jJ$ of Lemma \ref{J}. Then the
category $\chain{D}$ of chain complexes of left $D$-vector spaces
is equivalent to the category of left $(\lL/\jJ)$-comodules.
\end{corollary}
\begin{proof}
It follows from Theorems \ref{equivalencia-C} and \ref{Q_B-flat}.
\end{proof}

\subsection{The main example.}\label{main-example-1}

Here we will explain why Pareigis's example \cite{Pareigis:81},
even in the noncommutative case, always works. Thus, we will check
using the first statement of Theorem \ref{Q_B-flat} that the
cochain complex $Q_{\bullet}$ associated to the example of the
$R$-ring $A$ considered in \ref{main-example}, always satisfies
condition (2) of Theorem \ref{equivalencia-C}. In this way the
category $\chain{R}$ of chain complexes of left $R$-modules is
always equivalent to the category of left
$\bara{\lL(A)}$-comodules, where $\lL(A)$ is the left
$R$-bialgebroid described in Example \ref{main-example}.

Recall from Example \ref{main-example}, the $R$-ring $A=R\oplus Rt$ which is the  trivial generalized ring extension of $R$. Set $1_A=(1,0)$ and $\fk{t}=(0,t)$, so
we have $\fk{t}^2=0$. It is easily seen that the kernel of the multiplication of $A$, i.e. $K={\rm Ker}(A\tensor{R}A \to A)$ is free as a left and right
$R$-module with basis $\{\partial \fk{t},\fk{t}\partial \fk{t} \}$. In fact $K$ is a free $A$-module with rank one and basis $\partial\fk{t}$.
We summarize the properties of the cochain complex $Q_{\bullet}$, as follows.

\begin{proposition}\label{Q-ff}
The cochain complex $Q_{\bullet}$ associated to the trivial generalized ring
$A=R\oplus Rt$, fulfils the following properties:
\begin{enumerate}[(i)]
 \item For every $m \geq 2$, $Q_m$ is free as a left and right $R$-module with rank two, and its basis (on both sides) is given by the set
$\left\{\underset{}{} \fk{t}\partial\fk{t}\tensor{A}\cdots \tensor{A}\partial\fk{t},\, \partial\fk{t}\tensor{A}\cdots \tensor{A}\partial\fk{t} \right\}$.
\item $Q$ is a homotopically trivial complex.
\item $Q_C$ is faithfully flat module.
\end{enumerate}
\end{proposition}
\begin{proof}
$(i)$ This is  proved by induction on $m$. \\
$(ii)$ The homotopy is given by switching the dual basis. Let
$q_m=
\partial\fk{t}\tensor{A}\cdots \tensor{A}\partial\fk{t}$, ($(m-1)$-times) and $q_1=1_A$ be the generating element of
$Q_m$. Then we define a homotopy $h_m: Q_{m+1} \to Q_m$ by sending $q_{m+1} \mapsto \fk{t}q_{m}$ and $\fk{t}q_{m+1} \mapsto q_m$, $h_0$ is the first projection.\\
$(iii)$ The fact that $Q_C$ is flat follows from Theorem
\ref{Q_B-flat}, since we know that $Q_{\bullet}$ is exact and
splits either by Lemma \ref{mu-splits}, or by  item $(ii)$ and
\cite[Th\'eor\`eme 2.4.1]{Godement:1973}. Following the notations
of the proof of Theorem \ref{Q_B-flat}, we can easily see that
each right $T=\begin{pmatrix} R & R\\ 0 & R\end{pmatrix}$-module
$Q^{(m)}=\partial Q_{m-1} \oplus \fk{t}q_m R$ is isomorphic to
$eT$, where $e$ is the canonical idempotent of $T$. Henceforth,
the same argument of the last part of the proof of Theorem
\ref{Q_B-flat} serves to deduce that $Q_C$ is actually a
faithfully flat module.
\end{proof}

\begin{corollary}
Let $R$ be any $\Bbbk$-algebra and $A=R\oplus Rt$ its trivial
generalized extension. Consider the left $R$-bialgebroid
$\lL(A)$ described in Example \ref{main-example} and its
quotient $R$-coring $\bara{\lL(A)}$ by the left ideal
$\left\langle 1_{\lL(A)}(r\tensor{}1^o-1\tensor{}r^o)
\right\rangle$. Then the functor $Q\tensor{C}-$ establishes an
equivalence between the categories of chain complexes
of left $R$-modules and the category of left
$\bara{\lL(A)}$-comodules. In particular, if $R$ is a
commutative ring, then $Q\tensor{C}-$ establishes in fact a
monoidal equivalence.
\end{corollary}
\begin{proof}
The main claim is an immediate consequence of Proposition
\ref{Q-ff} and Theorem \ref{equivalencia-C}. In the last
statement, the functor in question  can be shown to be monoidal
using a similar proof of Lemma \ref{Q-monoidal}.
\end{proof}

\bigskip

\textbf{Acknowledgements} A. Ardizzoni would like to thank all the
members of the  Department of Algebra of the University of Granada
for the warm hospitality during his visit. L. El Kaoutit would like to thank all the
members of the  Department of Mathematics of University of Ferrara
for a warm hospitality during his visit.


\bigskip


\begin{thebibliography}{}


 \bibitem{Artin:1969}
M.~Artin, \emph{On Azumaya algebras and finite dimensional representations of rings}, J. Algebra, \textbf{11} (1969), 532--563.


\bibitem{Bourbaki:Chap. I-IV}
N.~Bourbaki, \emph{\'El\'ements de Math\'ematique. Alg\`ebre commutative. Chapitres 1 \`a 4}. Masson,  Paris 1985.


\bibitem{Brzezinski/Militaru:2002} T.~Brzezi\'nski and G.~Militaru, \emph{Bialgebroids, $\times_R$-Bialgebras and duality}.
J. Algebra, \textbf{251} (2002), 279--294.


\bibitem{Brzezinski/Wisbauer:2003}
T.~Brzezi\'{n}ski and R.~Wisbauer, \emph{Corings and comodules},
LMS, vol. 309, Cambridge University Press, 2003.


\bibitem{Bohm:arX08}
G.~B\"{o}hm, \emph{Hopf algebroids}.  Handbook of algebra.
Vol. \textbf{6}, 173--235.
Elsevier/North-Holland, Amsterdam, 2009.

\bibitem{Bruguieres:1994}
A.~Brugui\`{e}res, \emph{Théorie tannakienne non commutative}, Commun. in Algebra \textbf{22} (1994), 5817--5860.

\bibitem{Caenepeel/Groot/Vercruysse:2006}
S.~Caenepeel, E.~De Groot and J.~Vercruysse, \emph{Constructing
infinite comatrix corings from colimits}. Appl. Cat. Structures \textbf{14} no.5-6 (2006), 539--565.

\bibitem{Cartan:1958}
Henri Cartan, \emph{Homologie et cohomologie d'une alg\`ebre gradu\'ee.} S\'eminiare Henri Cartan, tome 11 n${}^o$. 2 (1958-1959), exp.
n${}^o$. 15, p. 1-20.

\bibitem{Deligne:1990}
P.~Deligne, \emph{Cat{\'e}gories tannakiennes}. In
\emph{The Grothendieck Festschrift} (P. Cartier et al., eds), Progr. math., 87, vol. II,
Birkh{\"a}user, Boston, MA. 1990, pp. 111--195.

\bibitem{Dold:1958}
A.~Dold, \emph{Homology of symmetric products and other functors of complexes}. Ann. Math. \textbf{68}, No. 1 (1958), 54--80.

\bibitem{El Kaoutit:2006}
L.~El Kaoutit, \emph{Corings over rings with local units}, Math. Nachr. \textbf{282} (2009), no.~5, 726--747.


\bibitem{El Kaoutit/Gomez:2003}
L.~El~Kaoutit and J.~G\'{o}mez-Torrecillas, \emph{Comatrix
corings: {G}alois coring, {D}escent theory, and a structure
theorem for cosemisimple corings}, Math. Z. \textbf{244} (2003),
887--906.

\bibitem{El Kaoutit/Gomez:2004}
L. El Kaoutit and J. G\'omez-Torrecillas, \emph{Infinite comatrix corings}. IMRN,
\textbf{39} (2004), 2017--2037.

\bibitem{ElKaoutit/Gomez/Lobillo:2004c}
L.~El~Kaoutit, J.~G\'{o}mez-Torrecillas, and F.~J. Lobillo,
\emph{Semisimple corings}, Algebra Colloq. \textbf{11} (2004),
no.~4, 427--442.

\bibitem{P.H. Hai:2008}
Ph\`{u}ng H\^{o} Hai, \emph{Tannaka-Krein duality for Hopf algebroids}, Israel J. Math. \textbf{167} (2008), 193--225.

\bibitem{Hovey:Model Catg}
M.~Hovey, \emph{Model categories}. Mathematical surveys and Monographs. Vol. \textbf{63}. AMS. 1999.




\bibitem{Gabriel:1962}
P.~Gabriel, \emph{Des cat\'{e}gories ab\'{e}liennes}, Bull. Soc.
Math. France \textbf{90} (1962), 323--448.

\bibitem{Godement:1973}
R. ~ Godement, \emph{Topologie alg\'ebrique et th\'eorie des faisceaux}. Troisi\`eme \'edition. Hermann, Paris 1973.


\bibitem{Gomez/Vercruysse:2005}
J.~G\'{o}mez-Torrecillas and J.~Vercruysse, \emph{Comatrix corings
and Galois comodules over firm rings}, Alg. Rep. Theory, \textbf{10} (2007), 271--306.

\bibitem{Goerss-Jardine}
Paul G.~Goerss, John F.~Jardine, \emph{Simplicial Homotopy Theory}. Progress in Mathematics, Vol. 174. Birkh\"auser,
Basel; Boston; Berlin. (1999).


\bibitem{Joyal/Street:1991}
A.~Joyal and R.~Street, \emph{An Introduction to Tannaka duality and quantum groups}. \emph{In ''Category Theory
Proceedings, Comom 1990 ''}, Lecture Note in Math. Vol. 1488, pp. 411--492, Springer-Verlag, Berlin 1991.





\bibitem{MacLane:Homology}
S. Mac Lane, \emph{Homology.} Grundlehren der math. Wissensch. \textbf{114}. Academic Press. Inc.,
Springer-Verlag, Berlin-G\"ottingen- Heidelberg, 1963.


\bibitem{McCrudden:2000}
Paddy McCrudden, \emph{Categories of Representations of coalgebroids}, Adv. Math. \textbf{154}, 299--332 (2000).


\bibitem{Pareigis:81}
B.~Pareigis, \emph{A non-commutative non-cocommutative Hopf algebra in ''nature''}, J. Algebra \textbf{70} (1981), 356--374.



\bibitem{Schauenburg:1998}
P.~ Schauenburg, \emph{Bialgebras over noncommutative rings and a structure theorem for Hopf bimodules}.
App. Catg. Struct. \textbf{6} (1998), 193--222.



\bibitem{Sweedler:1975b}
M. Sweedler, \emph{Groups of simple algebras}. I. H. E. S. Publ., \textbf{44} (1975), 79--189.



\bibitem{Takeuchi:1977b} M. Takeuchi, \emph{Groups of algebras over $A\tensor{}\bara{A}$}, J. Math. Soc. Japan \textbf{29} (1977), 459--492.

\bibitem{Tambara} D. Tambara, \emph{The coendomorphism bialgebra of an
algebra}, J. Fac. Sci. Univ. Tokyo Sect. IA Math. \textbf{37} (1990), no. 2,
425--456.

\bibitem{Wisbauer:2006}
R.~Wisbauer, \emph{On Galois comodules}, Commun. Algebra \textbf{34}
(2006), 2683--2711.


\end{thebibliography}
\end{document}